\documentclass[moor,nonblindrev,copyedit]{informs1}              % for a regular run
\usepackage{amsmath}
\usepackage{amssymb, amsmath, amsfonts, latexsym, verbatim}

\usepackage{natbib}
 \NatBibNumeric
 \def\bibfont{\small}%
 \bibpunct[, ]{[}{]}{,}{n}{}{,}%

\usepackage[colorlinks=true,breaklinks=true,bookmarks=true,urlcolor=blue,
     citecolor=blue,linkcolor=blue,bookmarksopen=false,draft=false]{hyperref}

\usepackage{amsmath}
\usepackage{amssymb, amsmath, amsfonts, latexsym}

\newcommand{\LA}{\mathcal{L}}

\catcode`\@=11

\newcommand{\A}{\mathcal{A}}
\newcommand{\B}{\mathcal{B}}

\newcommand{\M}{\mathcal{M}}

\newcommand{\NM}{\mathbb{N}}
\newcommand{\HL}{\mathcal{H}}
\newcommand{\R}{\mathbb{R}}

\newcommand{\F}{\mathcal{F}}

\newcommand{\SM}{\mathcal{S}}

\newcommand{\PR}{\mathcal{P}}
\newcommand{\defn}{\stackrel{def}{=}}
\newcommand{\ra}{\rightarrow}

\newcommand{\lra}{\longrightarrow}

\newcommand{\ep}{\epsilon}
\newcommand{\noi}{\noindent}

\def\EMAIL#1{\href{mailto:#1}{#1}}% When hyperref is used, otherwise outcomment
\def\URL#1{\href{#1}{#1}}         % When hyperref is used, otherwise outcomment

\TheoremsNumberedThrough     % Preferred (Theorem 1, Lemma 1, Theorem 2)
\EquationsNumberedThrough    % Default: (1), (2), ...
\MANUSCRIPTNO{MOR-2015-241} % When the article is logged in and DOI assigned to it,
\begin{document}
%%%%%%%%%%%%%%%%

% Outcomment only when entries are known. Otherwise leave as is and
%   default values will be used.
%\setcounter{page}{1}
%\VOLUME{00}%
%\NO{0}%
%\MONTH{Xxxxx}% (month or a similar seasonal id)
%\YEAR{0000}% e.g., 2005
%\FIRSTPAGE{000}%
%\LASTPAGE{000}%
%\SHORTYEAR{00}% shortened year (two-digit)
%\ISSUE{0000} %
%\LONGFIRSTPAGE{0001} %
%\DOI{10.1287/xxxx.0000.0000}%

% Author's names for the running heads
% Sample depending on the number of authors;
% \RUNAUTHOR{Jones}
% \RUNAUTHOR{Jones and Wilson}
% \RUNAUTHOR{Jones, Miller, and Wilson}
% \RUNAUTHOR{Jones et al.} % for four or more authors
% Enter authors following the given pattern:
\RUNAUTHOR{Basu and Ghosh}

% Title or shortened title suitable for running heads. Sample:
% \RUNTITLE{Bundling Information Goods of Decreasing Value}
% Enter the (shortened) title:
\RUNTITLE{Nonzero-sum Risk-sensitive Games with Countable States}

% Full title. Sample:
% \TITLE{Bundling Information Goods of Decreasing Value}
% Enter the full title:
\TITLE{Nonzero-sum Risk-sensitive Stochastic Games on a Countable State Space}

% Block of authors and their affiliations starts here:
% NOTE: Authors with same affiliation, if the order of authors allows,
%   should be entered in ONE field, separated by a comma.
%   \EMAIL field can be repeated if more than one author
\ARTICLEAUTHORS{%
\AUTHOR{Arnab Basu}
\AFF{Quantitative Methods and Information Systems Area, Indian
Institute of Management Bangalore, Bangalore- 560076, India, \EMAIL{arnab.basu@iimb.ernet.in}, \URL{http://www.iimb.ernet.in/user/112/basu-arnab}}
\AUTHOR{Mrinal K.\ Ghosh} \AFF{Department of Mathematics, Indian
Institute of Science, Bangalore - 560012, India,
\EMAIL{mkg@math.iisc.ernet.in},
\URL{http://math.iisc.ernet.in/~mkg}}
% Enter all authors
} % end of the block

\ABSTRACT{%
The infinite horizon risk-sensitive discounted-cost and ergodic-cost nonzero-sum stochastic games for controlled Markov chains with countably many states are analyzed. For the discounted-cost game, we prove the existence of Nash equilibrium strategies in the class of Markov strategies under fairly general conditions. %when the cost is independent of the controls of the players.
Under an additional geometric ergodicity condition and a small cost criterion, the existence of Nash equilibrium strategies in the class of stationary Markov strategies is proved for the ergodic-cost game. %via the analysis of the corresponding discounted-cost game.
% Enter your abstract
}%

% Sample
%\KEYWORDS{deterministic inventory theory; infinite linear programming duality;
%  existence of optimal policies; semi-Markov decision process; cyclic schedule}
%\MSCCLASS{Primary: 90B05; secondary: 90C40, 90C90}
%\ORMSCLASS{Primary: Inventory/production: deterministic multi-item;
%  secondary: dynamic programming/optimal control: deterministic
%  semi-Markov; programming: infinite dimensional}
%\HISTORY{Received November 20, 2003; revised March 8, 2004, and March 26, 2004.}

% Fill in data. If unknown, outcomment the field
 \MSCCLASS{Primary: 91A15, 91A25; secondary: 91A10, 91A50} \ORMSCLASS{Games/group
decisions: Noncooperative, Stochastic; Probability: Stochastic
model applications}
\HISTORY{Submitted September 16, 2015}

\maketitle
%%%%%%%%%%%%%%%%%%%%%%%%%%%%%%%%%%%%%%%%%%%%%%%%%%%%%%%%%%%%%%%%%%%%%%

\date{}

\section{Introduction} \label{intro}
We study risk-sensitive nonzero-sum stochastic games on the
infinite time horizon on a countable state space. Risk-sensitive
cost criterion plays an important role in many applications
including mathematical finance (see, e.g., Bilelecki and Pliska \cite{Biel1}, Nagai \cite{NH}).
In this criterion one investigates `exponential of integral' cost
which takes into account the attitude of the controller with
respect to risk. The study of this kind of cost criteria was first
initiated by Bellman \cite{Bell1}, p.\ 329 for the finite-state
space case. Howard and Matheson \cite{HoMa} did an in-depth
analysis for the first time in the finite-state space case where
each controlled chain is irreducible and aperiodic. Rothblum
\cite{Ro} extended it to the general non-irreducible finite-state
space case. In the past two decades, there has been a renewed
interest in this type of cost criteria as, when the `risk factor'
is strictly positive, i.e., in the risk-averse case, the use of
the exponential reduces the possibility of rare but devastating
large excursions of the state process. Though this criterion has
been studied extensively in the literature of Markov decision
processes (see, e.g., Borkar and Meyn \cite{BM1}, Cavazos-Cadena
and Fernandez-Gaucherand \cite{CCFG}, Di Masi and Stettner
\cite{DS1,DS1a,DS2}, Fleming and Hern\'{a}ndez-Hern\'{a}ndez
\cite{FH0}, Fleming and McEneaney \cite{FM1},
Hern\'{a}ndez-Hern\'{a}ndez and Marcus \cite{HM1}, Whittle
\cite{W0,W1}), the corresponding results on stochastic games seem
to be limited (see e.g., Basar \cite{Ba1}, El-Karoui and Hamadene
\cite{EKH}, Jacobson \cite{J1}, James et.\ al.\ \cite{JBE}, and,
Klompstra \cite{K1}). The general Linear-Exponential-Gaussian (LEG)
control problem for discrete time with perfect state observation
is treated in Jacobson \cite{J1} where an equivalence of this with
deterministic zero-sum quadratic-cost games was shown. Whereas this paper
addresses the undiscounted case, the corresponding discounted case was
addressed by Hansen and Sargent \cite{HS1}. This was further
extended to studying the Nash equilibrium for a two-person
discrete-time nonzero-sum game with quadratic-exponential cost
criteria in Klompstra \cite{K1}, the analogue of the
Linear-Exponential-Quadratic-Gaussian (LEQG) control problem
studied in Whittle \cite{W,W0}. The papers of Basar \cite{Ba1} and
El-Karoui et.\ al.\ \cite{EKH} deal with stochastic differential
games on the finite time horizon. In James et.\ al.\ \cite{JBE},
the finite-horizon risk-sensitive stochastic optimal control
problem for discrete-time nonlinear systems was studied and its
relation to a deterministic partially observed dynamic game was
established. To the best of our knowledge, the general case of
infinite-horizon risk-sensitive stochastic zero-sum games for both discounted
as well as ergodic cost criteria was first addressed by the current authors
for the differential games setup in Basu and Ghosh \cite{BaGh1} and
for the discrete-time countable state space case in Basu and Ghosh \cite{BaGh2}.
The case of nonzero-sum games in such setups still remains open.
In this work, we address this novel problem and analyze the generic nonzero-sum
case for the discrete-time countable state space setup using totally different
mathematical techniques for its solution as compared to the corresponding zero-sum case.
We study this problem with two players as the analysis can be routinely extended
to three or more player case without introducing any technical novelty but at the
cost of further notational complication.

We use the results of Balaji and Meyn \cite{BaM1} to extend the
work of Borkar and Meyn \cite{BM1} from one-controller case to
two-controller case in a fully competitive setup. In other words
we study nonzero-sum risk-sensitive stochastic games on the
infinite planning horizon with both discounted and ergodic cost
criteria. We would like to elucidate that both Balaji and Meyn
\cite{BaM1} as well as Borkar and Meyn \cite{BM1} have obtained
the desired results under a ``norm-like" or ``near-monotone"
condition on the running cost and a Lyapunov-type stability
condition. This ``norm-like" condition has been crucially used in
the proofs therein as to ensure that the `relative value
functions' are bounded from below. However, for our case, we have
performed the analysis without this ``norm-like" condition on the
running cost since it would not be suitable for our purpose as it
would invariably favour one of the competing players. This change
makes our analysis totally novel and substantially different from
those in the existing literature. Also, in most
of the existing literature in this domain see, e.g., Balaji and Meyn \cite{BaM1},
Borkar and Meyn \cite{BM1}, Cavazos-Cadena and
Fernandez-Gaucherand \cite{CCFG}, Di Masi and Stettner
\cite{DS1a}, and, Hern\'{a}ndez-Hern\'{a}ndez and Marcus
\cite{HM1}, the `risk factor' is assumed to be sufficiently small.
 We make an assumption, for the ergodic game only,  on the smallness of
 the cost function as in Basu and Ghosh \cite{BaGh1} which essentially
 implies that the `risk factor' cannot be too large.

%Using the approach of Borkar and Ghosh \cite{BG1}, we have
%used feedback relaxed strategies to study
%these stochastic games which maintain the strict
%competitive nature of the games at all times.
Under certain assumptions, we have established the existence of Nash equilibria
for both criteria. We obtain our results by studying the
corresponding Bellman equations. Note that if player I announces
that he is using a stationary/Markov strategy then for player II
the game problem reduces to a Markov decision problem (MDP). Then
by the results of Borkar and Meyn \cite{BM1}, Di Masi and Stettner \cite{DS1}
and Hern\'{a}ndez-Hern\'{a}ndez and Marcus \cite{HM1}, player II
has optimal stationary/Markov strategies. Such a strategy of
player II is called an `optimal response' of player II
corresponding to the announced strategy of player I. Optimal
responses of player I to the announced strategies of player II
are obtained analogously. Thus, for a given pair of strategies of
the two players, there exists a set of optimal responses. This
defines a point-to-set map. Any fixed point of this map is clearly
a Nash equilibrium. In this paper, we establish the existence of
such a fixed point thereby establishing the existence of Nash
equilibria for relevant cases. For the sake of notational
simplicity, we consider two-player games only. All our results
extends to multi-player games in a routine manner.

The rest of our paper is structured as follows. Section \ref{probdescpr} deals
with the description of the problem. The discounted cost criterion
is studied in Section \ref{nonzero-sum-disc}. Section
\ref{nonzero-sum-erg} deals with the ergodic cost criterion. %We also provide
%an illustrative example in this section.
We conclude our paper in Section \ref{conclude} with a summary and possible future directions of work.

\section{Problem Description} \label{probdescpr}

A two-person stochastic game is determined by six objects
$(X,U,V,r_1,r_2,q)$ where $X \defn \{0,1,2,\ldots\}$ is a
countable state space; $U$ and $V$ are action spaces of players
I and II, resp., assumed to be compact metric spaces; $r_1 (\mbox{resp.}\ r_2) : X
\times U \times V \mapsto \R$ is the one-stage cost function
for player I (\mbox{resp.}\ II) assumed to be bounded and jointly continuous in
$(u,v) \in U \times V$ for each $k \in X$. Let $\PR(X)$ be the space of
probability measures on $X$ endowed with the Prohorov topology
(see, e.g., Borkar \cite{B1}). Let $q : X \times U \times V \mapsto \PR(X)$
be the transition stochastic kernel which is assumed to be jointly continuous in $(u,v) \in U \times V$
in the topology of weak convergence for each $k \in X$.
The game is played is as follows: At each stage (time instant)
players observe the current state $k \in X$ of the system and then
players I and II independently choose actions $u \in U$, $v \in
V$, resp. As a result of this, two things happen:
\begin{itemize}
\item[(i)] the player I (resp.\ II)  pays an immediate cost
$r_1(k,u,v) (\mbox{resp.}\ r_2(k,u,v))$,
\item[(ii)] the system moves to a new state $k' \in
X$ with distribution $q(\cdot|k,u,v)$.
\end{itemize}
The whole process then repeats from the new state $k'$. Cost
accumulates throughout the course of the game. The planning
horizon or total number of stages is infinite, and each player
wants to minimize his infinite-horizon multiplicative expected
cost to be described shortly. At each stage the players choose
their actions independently on the basis of past information. The
available information for decision making at time $t \in \NM_0 \defn
\{0,1,2\ldots\}$ is given by the history of the process up to that
time \[h_t \defn (k_0, (u_0,v_0), k_1, (u_1,v_1), \ldots,
(u_{t-1},v_{t-1}), k_t) \in H_t,\] where $H_0 = X,\ H_t = H_{t-1}
\times (U \times V \times X), \ldots, H_{\infty} = (U \times V
\times X)^{\infty}$ are the history spaces. A strategy for player
I is a sequence $\mu \defn \{\mu_t : H_t \mapsto \PR(U)\}_{t \in
\NM}$ of stochastic kernels.  The set of all such strategies for
player I is denoted by $\Pi_1$.
%up to and including time $T$ is given by $\Pi_1^T \defn \{\mu_t :
%H_t \ra \PR(U)\}_{0 \leq t \leq T}$. Similarly, we define $\hat{\Pi}_1^T \defn \{\mu_t :
%H_t \ra \PR(U)\}_{t \geq T}$ to be the set of strategies starting at time $T$.
Given any $t,T \in \NM_0$ with $T \geq t+1$, the set of strategies
$\{\mu_s : H_s \mapsto \PR(U)\}_{t \leq s \leq T}$ played from time $t$ to time $T$ is denoted by $\Pi_1^{t,T}$.
Hence, $\Pi_1^{0,\infty} \equiv \Pi_1$.

\begin{definition} \label{marstatstr}
A strategy $\mu$ for player
I is called a Markov strategy if
\[\mu_t(h_{t-1},u,v,k)(\cdot) = \mu_t(h'_{t-1},u',v',k)(\cdot)\]
for all $h_{t-1},h'_{t-1} \in H_{t-1}$, $u,u' \in U$, $v,v' \in
V$, $k \in X$, $t \in \NM_0$. Thus a Markov strategy for player I
can be identified with a sequence of measurable maps, denoted by
$\mu \equiv \{\mu_t : X \mapsto \PR(U)\}_{t \in \NM}$. A Markov strategy
$\{\mu_t\}$ is called stationary Markov if $\mu_t = \mu : X \mapsto
\PR(U)$ for all $t$. A stationary Markov strategy is called pure
or deterministic if $\mu : X \mapsto U$.
\end{definition}

Let $\M_1, \SM_1, D_1$ denote
the set of Markov, stationary Markov and deterministic strategies
strategies for player I. The strategy sets $\Pi_2, \Pi_2^{t,T}, %\Pi_2^{T},
\M_2,\SM_2,D_2$ for player II are defined analogously. The spaces
$\M_i,\SM_i,\ i=1,2$ are endowed with the product topologies derived
from the Prohorov topology on the underlying spaces $\PR(U)$ (resp.\ $\PR(V)$)
respectively. Since $U,V$ are compact metric spaces, it follows that $\M_i,\SM_i,\
i=1,2$ are also compact metric spaces. Note that going forward we sometimes
use barbarism of notation and denote a stationary strategy $\{\mu,\mu,\ldots\}$
as $\{\mu\}$ or only $\mu$ when the context is clear.

\begin{comment}
We assume that the following condition holds throughout this paper: \\
\noi \textbf{(A0)} The one-stage cost function $r_i,\ i=1,2$  is independent
of $(u,v) \in U \times V$.

For $\xi \in \PR(U), \chi \in \PR(V)$, we use the notation
%\begin{equation} \label{avgruncost}
%r_i(k,\xi,\chi) \defn \int_V \int_U r_i(k,u,v) \xi(du)
%\chi(dv),\ i=1,2,
%\end{equation}
%and
\begin{equation} \label{avgtranker}
\hat{q}(k'|k,\xi,\chi) \defn \int_V \int_U q(k'|k,u,v) \xi(du)
\chi(dv),\ i=1,2.
\end{equation}
\end{comment}

Given an initial distribution $\pi_0 \in \PR(X)$ and a pair of
strategies $(\mu,\nu) \in \Pi_1 \times \Pi_2$, the corresponding
state and action processes $\{X_t\}, \{U_t\}, \{V_t\}$ are
stochastic processes defined on the canonical space $(H_{\infty},
\B(H_{\infty}), P^{\mu,\nu}_{\pi_0})$ ($\B(H_{\infty}) =$ Borel
$\sigma$-field on $H_{\infty}$) via the projections
$X_t(h_{\infty}) = k_t, U_t(h_{\infty}) = u_t, V_t(h_{\infty}) =
v_t$, where $P^{\mu,\nu}_{\pi_0}$ is uniquely determined by
$\mu,\nu$ and $\pi_0$ by Ionescu Tulcea's Theorem
(see, e.g., Proposition 7.28 of Bertsekas and Shreve \cite{BS1}).
When $\pi_0 = \delta_k,\ k \in X$, we simply write
$P_k^{\mu,\nu}$.

Let $(\mu,\nu) \in \Pi_1 \times \Pi_2$ and let $(X^{(t)}_s,U^{(t)}_s,V^{(t)}_s: s\geq t)$ be
the corresponding process starting from $X^{(t)}_t = k \in X$ for given $t \geq 0$. We omit
the superscript when $t=0$ for notational convenience.

\begin{definition} \label{costcrit}
The risk-sensitive
discounted cost for player I ($i=1$) and II ($i=2$) is defined by
\begin{equation} \label{riskdiscpayoff}
\rho_i^{\mu,\nu}(\theta,(k,t)) = \frac{1}{\theta} \ln E_{k,t}^{\mu,\nu}\left[e^{\theta \sum_{s=t}^{\infty}\alpha^{s-t} r_i(X^{(t)}_s,U^{(t)}_s,V^{(t)}_s)}\right],
\end{equation}
where $\theta \in (0,\Theta],\ \Theta > 0$ is the `risk-sensitive'
parameter, $\alpha \in [0,1]$ is the `discount factor' and
$E_{k,t}^{\mu,\nu}$ denotes the expectation with respect to
$P^{\mu,\nu}_{k,t}$. When $t=0$, we omit this subscript `t' for notational convenience.
The corresponding risk-sensitive ergodic cost for
player $i,\ i=1,2$ is defined by
\begin{equation} \label{riskergpayoff}
\beta_i^{\mu,\nu}(\theta,k) = {\lim\sup}_{T \ra \infty}\frac{1}{\theta T} \ln E_k^{\mu,\nu}\left[e^{\theta \sum_{t=0}^{T-1} r_i(X_t,U_t,V_s)}\right].
\end{equation}
\end{definition}

\begin{definition} \label{equil}
Given a $(\theta,k) \in (0,\Theta] \times X$, a pair of strategies
$(\mu^*,\nu^*) \in \Pi_1 \times \Pi_2$ is called a Nash
equilibrium (for the cost criteria above) if
\begin{eqnarray} \label{Nashdefn}
&& L_1^{\mu^*,\nu^*}(\theta,k) \leq L_1^{\mu,\nu^*}(\theta,k),\ \mu \in \Pi_1, \nonumber \\
&& \mbox{and} \nonumber \\
&& L_2^{\mu^*,\nu^*}(\theta,k) \leq L_2^{\mu^*,\nu}(\theta,k),\
\nu \in \Pi_2,
\end{eqnarray}
where for $i=1,2$, $L_i = \rho_i$ or $L_i = \beta_i$ as the case
may be.
\end{definition}

For a given $\theta \in (0,\Theta]$, a pair
$(\mu^*(\theta,\cdot),\nu^*(\theta,\cdot))$ of stationary/Markov
strategies (depending on $\theta$) is said to be Nash equilibrium
strategies if these measurable maps constitute a Nash equilibrium
for any initial $k \in X$.

%%%%%%%%%%%%% TO BE REMOVED %%%%%%%%%%%%%%%%%%%%%%%%
%We now state a result which shall be used in Section
%\ref{nonzero-sum-erg}. For its proof, refer to Proposition II.4.2 of
%\cite{DE1}.

%\begin{proposition} \label{dupuis-ellis-theorem}
%Let $f \in B(X)$ and let $\xi \in \PR(X)$. Then
%\begin{equation*}
%\ln \sum_{j \in X} e^{f(j)} \xi(j) = \sup_{\chi \in
%\PR(X)}\left\{\sum_{j \in X} f(j) \chi(j) - I[\chi||\xi]\right\},
%\end{equation*}
%where, for any two probability measures $p,q$,
%\begin{equation*}
%I[p||q] \equiv  \left\{ \begin{array}{ll} \sum_{j \in X}p(j)\log\frac{p(j)}{q(j)} & \mbox{if $p \ll q$} \\ +\infty & \mbox{otherwise} \end{array} \right.
%\end{equation*}
%denotes their relative entropy. Furthermore, the supremum is attained at a
%unique probability measure $\chi^* \in \PR(X)$ given by
%\begin{equation*}
%\chi^*(k) = \frac{e^{f(k)}}{\sum_{j \in X} e^{f(j)} \xi(j)}
%\xi(k),\ k \in X.
%\end{equation*}
%\end{proposition}
%%%%%%%%%%%%%%%%%%%%%%%%%%%%%%%%%%%%%%%%%%%%%%%%%%%%%

%Under \textbf{(A0)}
We shall first establish the existence of Nash equilibria in the class of Markov strategies for cost criterion (\ref{riskdiscpayoff}).
 Under additional ergodicity and smallness of cost conditions, we establish the existence of Nash equilibria in the class of stationary Markov strategies for cost criterion (\ref{riskergpayoff}).

Given a topological space $Y$, we denote by $B(Y)$ and $C_b(Y)$
the Banach spaces of bounded measurable and bounded continuous
functions on $Y$ respectively, each equipped with the sup-norm
metric.

The following result established in Section II.3, Chapter II (pp.\ 25-30) of Borkar \cite{B0} plays a crucial role
in Sections \ref{nonzero-sum-disc} and \ref{nonzero-sum-erg}. Given any set of stochastic process
$\{Y^{(1)}_t\},\{Y^{(2)}_t\},\ldots$ on $(H_{\infty},\B(H_{\infty}))$ driven by $(\mu,\nu) \in \Pi_1 \times \Pi_2$
we denote their joint  law by $\LA^{\mu,\nu}\{(Y_t^{(1)},Y_t^{(2)},\ldots), t\geq 0\}$.

\begin{proposition} \label{contlawstrat}
For a fixed initial distribution $\pi_0 \in \PR(X)$, and given
$(\mu,\nu) \in \M_1 \times \M_2$, the map \[\M_1 \times \M_2 \ni
(\mu,\nu) \mapsto \LA^{\mu,\nu}\{(X_t,U_t,V_t), t \geq 0\}\]
is jointly continuous in $(\mu,\nu)$, where $\{X_t, t \geq 0\}$ denotes the state
process with initial law $\pi_0$ and $\{U_t, t \geq 0\}$, $\{V_t,t \geq 0\}$
are the corresponding control processes.
\end{proposition}

%\noi \textbf{Proof:}

\section{Discounted Game} \label{nonzero-sum-disc}

In this section, we study the cost criterion (\ref{riskdiscpayoff}).
To this end, we first consider the risk-sensitive exponential cost
criterion for player I ($i=1$) and II ($i=2$) given by:
\begin{equation} \label{riskexppayoff}
\zeta^{\mu,\nu}_i(\theta,(k,t)) \defn
e^{\theta\rho^{\mu,\nu}_i(\theta,(k,t))} =
E_{k,t}^{\mu,\nu}\left[e^{\theta
\sum_{s=t}^{\infty} \alpha^{s-t} r_i(X_s^{(t)},U_s^{(t)},V_s^{(t)})}\right],
\end{equation}
where $\{X_s^{(t)}:\ s \geq t\}$ is the state process and $\{(U_s^{(t)},V_s^{(t)}):\ s \geq t\}$ are the control processes under $(\mu,\nu) \in \Pi_1^{t,\infty} \times \Pi_2^{t,\infty}$ starting from $k \in X$.

Given strategies $(\mu,\nu) \in \Pi^{t,\infty}_1 \times \Pi^{t,\infty}_2$, consider the
following Bellman equations for the
exponential cost (\ref{riskexppayoff}) for players I and II (resp.):
\begin{eqnarray} \label{HJIexpupper}
\phi_1(\theta,(k,t))= \inf_{\xi \in \PR(U)}
\left[\int_U \int_V e^{\theta r_1(k,u,v)} \left(\sum_{j \in X}
\phi_1(\theta\alpha,(j,t+1)) q(j|k,u,v)\right) \xi(du)\nu_t[h_t](dv)\right],
\end{eqnarray}
with
\begin{equation} \label{HJIexpbnd}
\lim_{\theta \ra 0} \phi_1(\theta,(k,t)) = 1,\ \forall k \in X,\ \forall t,
\end{equation}
and
\begin{eqnarray} \label{HJIexplower}
\phi_2(\theta,(k,t)) =  \inf_{\chi \in \PR(V)}
\left[\int_U \int_V e^{\theta r_2(k,u,v)} \left(\sum_{j \in X}
\phi_2(\theta\alpha,(j,t+1)) q(j|k,u,v)\right)\mu_t[h_t](du)\chi(dv))\right],
\end{eqnarray}
with a boundary condition analogous to (\ref{HJIexpbnd}).

Note that under $(\mu,\nu) \in \Pi^{t,\infty}_1 \times \Pi^{t,\infty}_2$, the corresponding chain $\{X_t\}$ is inhomogeneous. Hence, we consider the transformed chain $\{\tilde{X}_t\}$ on $X \times \NM_0$ with the transition kernel $\tilde{q}$ defined by
\begin{equation} \label{transker}
\tilde{q}((k,t')|(j,t),u,v) = \hat{q}(k|j,u,v) \delta_{t',t+1},\ u \in U,\ v \in V,
\end{equation}
where $\delta_{s',s} =1$ if $s'=s$ and $0$ otherwise. Now, replacing $q$ by $\tilde{q}$ and using the arguments used in the proof of Proposition 3.1 of Di Masi and Stettner \cite{DS1}, we get the following result. We omit the routine details.

%We now state and prove the following proposition.

\begin{proposition} \label{dynprogsolexp}
Given $(\mu,\nu) \in \Pi^{t,\infty}_1 \times \Pi^{t,\infty}_2$, there exist unique bounded solutions $\hat{\phi}_1[\nu],\hat{\phi}_2[\mu]$ to
(\ref{HJIexpupper}) and (\ref{HJIexplower}) (resp.) satisfying the
boundary condition (\ref{HJIexpbnd}) such that
\begin{equation} \label{stochrep01}
\hat{\phi}_1[\nu](\theta,(k,t)) = \inf_{\tilde{\mu} \in
\Pi_1^{t,\infty}} \zeta^{\tilde{\mu},\nu}_1(\theta,(k,t)),
\end{equation}
and
\begin{equation} \label{stochrep02}
\hat{\phi}_2[\mu](\theta,(k,t)) = \inf_{\tilde{\nu} \in
\Pi_2^{t,\infty}} \zeta^{\mu,\tilde{\nu}}_2(\theta,(k,t)).
\end{equation}
\end{proposition}

Moreover, by Remark 3.2 of Di Masi and Stettner \cite{DS1}, given $(\mu,\nu) \in \Pi^{t,\infty}_1 \times \Pi^{t,\infty}_2$ and $\theta \in (0,\Theta]$, the minimizing strategies
$\{\mu_t^*[\nu]\},\ \{\nu_t^*[\mu]\} $ in
(\ref{stochrep01}) and (\ref{stochrep02}) (resp.) are given by
\begin{equation} \label{minseq01}
\mu^*_t[\nu] = \hat{\mu}[\nu](\theta \alpha^t,(X_t,t)) \in \M_1
\end{equation}
and
\begin{equation} \label{minseq02}
\nu^*_t[\mu] = \hat{\nu}[\mu](\theta \alpha^t,(X_t,t)) \in \M_2
\end{equation}
%By the continuity of $q$ and the compactness of $\PR(U)$ and $\PR(V)$, there exists measurable selectors (see Bene\v{s} \cite{Be1})
where $(\hat{\mu}[\nu],\hat{\nu}[\mu]) : (0,\Theta] \times (X \times \NM_0) \mapsto \PR(U) \times \PR(V)$ are measurable (minimizing) selectors (see Bene\v{s} \cite{Be1}) such that
\begin{eqnarray} \label{discminsel1}
\inf_{\xi \in \PR(U)} \left[\int_U \int_V e^{\theta r_1(k,u,v)} \left(\sum_{j \in X} \hat{\phi}_1[\nu](\theta\alpha,(j,t+1)) q(j|k,u,v)\right) \xi(du)\nu_t[h_t](dv)\right] \nonumber \\
= \int_U \int_V e^{\theta r_1(k,u,v)} \left(\sum_{j \in X} \hat{\phi}_1[\nu](\theta\alpha,(j,t+1))q(j|k,u,v)\right)\hat{\mu}[\nu](\theta,(k,t))(du)\nu_t[h_t](dv),
\end{eqnarray}
and
\begin{eqnarray} \label{discminsel2}
\inf_{\chi \in \PR(V)} \left[\int_U \int_V e^{\theta r_2(k,u,v)} \left(\sum_{j \in X} \hat{\phi}_2[\mu](\theta\alpha,(j,t+1)) q(j|k,u,v)\right)\mu_t[h_t](du)\chi(dv)\right] \nonumber \\
= \int_U \int_V e^{\theta r_2(k,u,v)} \left(\sum_{j \in X} \hat{\phi}_2[\mu](\theta\alpha,(j,t+1))q(j|k,u,v)\right)\mu_t[h_t](du)\hat{\nu}[\mu](\theta,(k,t))(dv).
\end{eqnarray}

Thus, $\mu^*_t[\nu] \in \M_1$ (resp.\ $\nu^*_t[\mu] \in \M_2$) is an optimal response of
player I (resp.\ player II) corresponding to $\nu \in \Pi^{t,\infty}_2$
(resp.\ $\mu \in \Pi^{t,\infty}_1$). Hence, without loss of generality, we can effectively
consider the strategy pair $(\mu,\nu) \in \M_1 \times \M_2$ for further analysis of this game.

Considering the minimizing selectors
$(\hat{\mu}[\nu],\hat{\nu}[\mu])$ in (\ref{discminsel1}) and (\ref{discminsel2}) (resp.), we now define the
following point-to-set maps $\HL_i : \M_j \mapsto 2^{\M_i},\
i,j=1,2,\ i\neq j$ as follows:
\begin{eqnarray} \label{optresponsedisc01}
\HL_1(\nu) \defn \{\{\mu^*_t[\nu]\} : \mu^*_t[\nu]\ \mbox{satisfies}\ (\ref{minseq01})\},
\end{eqnarray}
and
\begin{eqnarray} \label{optresponsedisc02}
\HL_2(\mu) \defn \{\{\nu^*_t[\mu]\} : \nu^*_t[\mu]\ \mbox{satisfies}\ (\ref{minseq02})\}.
\end{eqnarray}

Define the map $\HL \equiv \HL_1 \otimes \HL_2 : \M_1 \times \M_2 \mapsto 2^{\M_1 \times \M_2}$ as
\begin{eqnarray} \label{optresponsedisc}
\HL(\mu,\nu) \defn \{(\{\mu^*_t[\nu]\},\{\nu^*_t[\mu]\}) : \mu^*_t[\nu]\ \mbox{satisfies}\ (\ref{minseq01})\ \mbox{and}\ \nu^*_t[\mu]\ \mbox{satisfies}\ (\ref{minseq02})\}.
\end{eqnarray}

%Now we state and prove an important result.

%\begin{proposition} \label{Nashlem01}
%Given $\theta \in (0,\Theta]$, the map $\HL$ has a fixed point $(\tilde{\mu},\tilde{\nu}) \in \M_1 \times \M_2$.
%\end{proposition}

%\noi \textbf{Proof:} \hfill $\Box$

We now prove the existence of Nash equilibria in the class of Markov strategies for the exponential cost criterion (\ref{riskexppayoff}).

\begin{theorem} \label{existexpvalue}
Given $\theta \in (0,\Theta]$, there exists a pair of Nash equilibrium strategy in $\M_1 \times \M_2$ for the game corresponding to the cost criterion (\ref{riskexppayoff}).
\end{theorem}

\noi \textbf{Proof:} Given $\theta \in (0,\Theta],\ (\mu,\nu) \in
\M_1 \times \M_2$, both $\HL_1(\nu)$ and $\HL_2(\mu)$ are
non-empty and convex implying $\HL(\mu,\nu)$ is also non-empty
and convex. %(as \textbf{(A0)} holds).
Let $\{(\{\mu^*_{t,k}[\nu]\}_t,\{\nu^*_{t,k}[\mu]\}_t)\}_{k \in \NM}$ be a sequence of distinct points of $\HL(\mu,\nu)$
converging to $(\{\mu_{t,\infty}[\nu]\}_t,\{\nu_{t,\infty}[\mu]\}_t)$, i.e., $(\{\mu_{t,\infty}[\nu]\}_t,\{\nu_{t,\infty}[\mu]\}_t)$
is a limit point of $\HL(\mu,\nu)$. Since, for each $k \in \NM$, $\mu^*_{t,k}[\nu] = \hat{\mu}_k[\nu](\theta \alpha^t,(X_t,t))$
and $ \nu^*_{t,k}[\mu] = \hat{\nu}_k[\mu](\theta \alpha^t,(X_t,t))$ satisfy (\ref{discminsel1}) and (\ref{discminsel2}) resp.,
it follows (by linearity) that so do $\mu_{t,\infty}[\nu]$ and $\nu_{t,\infty}[\mu]$ resp. implying that they satisfy (\ref{minseq01}) and (\ref{minseq02}) resp. The closure property of $\HL(\mu,\nu)$ thus follows.
% from the joint continuity of $q(\cdot|\cdot,u,v)$ under $(\mu,\nu)$.
Hence, $\HL$ is a map with non-empty, closed and convex values. Now, by
Proposition \ref{contlawstrat}, $\LA^{\mu,\nu}\{X_t, t \geq 0\}$
is jointly continuous in $(\mu,\nu) \in \M_1 \times \M_2$. Note that $\hat{\phi}_1[\nu]$ as obtained in
(\ref{stochrep01}) has a minimizing Markov strategy
(\ref{minseq01}) and $\hat{\phi}_2[\mu]$ as obtained in
(\ref{stochrep02}) has a minimizing Markov strategy
(\ref{minseq02}), i.e., the minimizations in (\ref{stochrep01}) and (\ref{stochrep02})
can be effectively considered over $\M_1$ and $\M_2$ (resp.) which are compact sets. Also, $r_i,\ i=1,2$ are bounded functions.
Hence, %by Proposition \ref{contlawstrat} and Skorohod's theorem (see Theorem 2.2.2 of Borkar \cite{B1}),
by Berge Maximum Theorem (see Theorem 17.31 of Aliprantis and Border \cite{AliBor1}), $\hat{\phi}_1[\nu]$ and
$\hat{\phi}_2[\mu]$ are continuous in $\nu$ and $\mu$
respectively. Now consider a sequence $\{(\mu_k,\nu_k)\}_{k \in
\mathbb{N}}$ in $\M_1 \times \M_2$ converging to
$(\mu_\infty,\nu_\infty) \in \M_1 \times \M_2$. Let
$\HL(\mu_k,\nu_k) \ni (\tilde{\mu}_k,\tilde{\nu}_k) \ra
(\tilde{\mu}_\infty,\tilde{\nu}_\infty) \in \M_1 \times \M_2$.
From the continuity of
$q(\cdot|\cdot,u,v),\hat{\phi}_1[\nu]$ and
$\hat{\phi}_2[\mu]$ in $\mu,\nu$, it follows that
$(\tilde{\mu}_\infty,\tilde{\nu}_\infty) \in
\HL(\mu_\infty,\nu_\infty)$. Hence the map $\HL$ is upper
semi-continuous. Then, by Theorem 1 of Fan \cite{F1}, there exists
a fixed point $(\mu^*,\nu^*) :\NM_0 \times X \mapsto \PR(U) \times
\PR(V)$ of the map $\HL$ defined as
\begin{equation} \label{expdisccostfp}
(\mu^*(t,k),\nu^*(t,k)) \defn (\hat{\mu}[\nu^*](\theta \alpha^t,(k,t)),\hat{\nu}[\mu^*](\theta\alpha^t,(k,t)))
\end{equation}
where $\hat{\mu}[\nu^*]$ (resp.\ $\hat{\nu}[\mu^*]$) is any minimizing selector in (\ref{discminsel1}) (resp. (\ref{discminsel2})).
Then $(\mu^*,\nu^*) \in \M_1 \times \M_2$ is a
Nash equilibrium strategy for the cost criterion (\ref{riskexppayoff}), i.e.,
\[\hat{\phi}_1[\nu^*](\theta,(k,t)) = \zeta^{\mu^*,\nu^*}_1(\theta,(k,t)) \leq \zeta^{\mu,\nu^*}_1(\theta,k),\ \forall \mu \in \Pi_1^{t,\infty}, \]
and
\[\hat{\phi}_2[\mu^*](\theta,(k,t)) = \zeta^{\mu^*,\nu^*}_2(\theta,(k,t)) \leq \zeta^{\mu^*,\nu}_2(\theta,k),\ \forall \nu \in \Pi_2^{t,\infty}. \]
\hfill $\Box$ % The result follows directly from the discussions above. \hfill $\Box$

Now we prove the existence of Nash equilibria in the class of Markov strategies for the cost criterion (\ref{riskdiscpayoff}).

Given definition (\ref{riskexppayoff}), the Bellman equations for the discounted cost (\ref{riskdiscpayoff}) are:
\begin{equation} \label{HJIdiscupper}
e^{\theta\psi_1(\theta,(k,t))} = \inf_{\xi \in \PR(U)} \left[\int_U \int_V e^{\theta r_1(k,u,v)} \left(\sum_{j \in X} e^{\theta \alpha \psi_1(\theta\alpha,(j,t+1))} q(j|k,u,v)\right) \xi(du)\nu_t[h_t](dv)\right],
\end{equation}
for player $1$ with the boundary condition
\begin{equation} \label{HJIdiscupperbnd}
\lim_{\theta \ra 0} \psi_1(\theta,(k,t)) = \inf_{\mu \in \Pi_1^{t,\infty}} E_{k,t}^{\mu,\nu}\left[\sum_{s=t}^{\infty}\alpha^{s-t} r_1(X^{(t)}_s,U^{(t)}_s,V^{(t)}_s)\right],\ \forall k \in X,\ \forall t,
\end{equation}
and
\begin{equation} \label{HJIdisclower}
e^{\theta\psi_2(\theta,(k,t))} = \inf_{\chi \in \PR(V)} \left[\int_U \int_V e^{\theta r_2(k,u,v)} \left(\sum_{j \in X} e^{\theta \alpha \psi_2(\theta\alpha,(j,t+1))} q(j|k,u,v)\right)\mu_t[h_t](du)\chi(dv)\right],
\end{equation}
for player $2$ with the boundary condition
\begin{equation} \label{HJIdisclowerbnd}
\lim_{\theta \ra 0} \psi_2(\theta,(k,t)) = \inf_{\nu \in \Pi_2^{t,\infty}} E_{k,t}^{\mu,\nu}\left[\sum_{s=t}^{\infty}\alpha^{s-t} r_2(X^{(t)}_s,U^{(t)}_s,V^{(t)}_s)\right],\ \forall k \in X,\ \forall t.
\end{equation}
%where $(X^{\mu,\nu}_s,U^{\mu,\nu}_s,V^{\mu,\nu}_s: s \geq t)$ is the process corresponding to $(\mu,\nu) \in \Pi_1^{t,\infty} \times \Pi_2^{t,\infty}$.

\begin{theorem} \label{existdisc}
Given $\theta \in (0,\Theta]$, consider the fixed point Markov strategies $(\mu^*,\nu^*)$ of $\HL$ as in (\ref{expdisccostfp}) in Theorem \ref{existexpvalue}. Then $(\mu^*,\nu^*)$ is also a pair of Nash equilibrium strategy for the
cost criterion (\ref{riskdiscpayoff}), i.e.,
\[\hat{\psi}_1[\nu^*](\theta,(k,t)) = \rho^{\mu^*,\nu^*}_1(\theta,(k,t)) \leq \rho^{\mu,\nu^*}_1(\theta,(k,t)),\ \forall \mu \in \Pi_1^{t,\infty}, \]
and
\[\hat{\psi}_2[\mu^*](\theta,(k,t)) = \rho^{\mu^*,\nu^*}_2(\theta,(k,t)) \leq \rho^{\mu^*,\nu}_2(\theta,(k,t)),\ \forall \nu \in \Pi_2^{t,\infty}, \]
\end{theorem}
where  $\hat{\psi}_1[\nu^*]$ and $\hat{\psi}_2[\mu^*]$ are unique bounded solutions to
(\ref{HJIdiscupper}) and (\ref{HJIdisclower}) (resp.) satisfying the
boundary conditions (\ref{HJIdiscupperbnd}) and (\ref{HJIdisclowerbnd}) respectively.

\noi \textbf{Proof:} The result follows directly from (\ref{riskexppayoff}), Proposition \ref{dynprogsolexp}, Theorem \ref{existexpvalue} and the fact that log is an increasing function. \hfill $\Box$

\section{Ergodic Game} \label{nonzero-sum-erg}
In this section, we study the cost criterion (\ref{riskergpayoff}). To this
end, we make the following assumption:\\
\noi \textbf{(A1)} The process $\{X_t\}_{t \in \NM_0}$ is an irreducible, aperiodic Markov chain under any pair of stationary Markov strategies $(\mu,\nu) \in \SM_1 \times \SM_2$.

We also assume the following condition which guarantees uniform ergodicity of controlled Markov processes
and is used to study additive average cost problems (see, e.g., Hern\'{a}ndez-Lerma \cite{HL1}, Section 3.3 or an equivalent assumption (A1) in Di Masi and Stettner \cite{DS1} or (3.2) in Hern\'{a}ndez-Hern\'{a}ndez and Marcus \cite{HM1}):\\
%(see, e.g., \cite{HL1}, Section 3.3):\\
\noi \textbf{(A2)} There exists a number $0 < \delta < 1$ such that for all $i,j \in X$, $u,u' \in U$ and $v,v' \in V$,
\[||q(\cdot|i,u,v) - q(\cdot|j,u',v')||_{\mbox{TV}} \leq 2 \delta,\] where $||\cdot||_{\mbox{TV}}$ denotes the total variation norm (see, e.g., Hern\'{a}ndez-Lerma \cite{HL1}, Appendix B).  %This condition is equivalent to the following (see, e.g., assumption (A1) in \cite{DS1}):\\
%\noi \textbf{(A2)} There exists a $0 < \delta < 1$ such that for all $i,j,k \in X$, $\xi,\xi' \in \PR(U)$ and $\chi,\chi' \in \PR(V)$,
%\[\hat{q}(k|i,\xi,\chi) - \hat{q}(k|j,\xi',\chi') \leq \delta.\]

In our further analysis of this problem, the above assumptions always hold.

%We also assume the following Doeblin-type condition (see, e.g., \cite{MT1}) to hold throughout this section which guarantees the stability of all stationary Markov strategies $(\mu,\nu) \in \SM_1 \times \SM_2$:\\
%\noi \textbf{(C0)} There exists a finite set $K \subset X$, an
%integer $m \geq 1$ and a number $R > 0$ such that under any
%$(\mu,\nu)  \in \SM_1 \times \SM_2$, $\hat{q}$ satisfies the
%following condition:
%\[\sum_{j\in K}\hat{q}^m(j|k,\mu,\nu) > R,\ \forall k \in X,\]
%where $\hat{q}^m(\cdot|k,\mu,\nu)$ denotes the $m$-step transition function under $(\mu,\nu)$ starting from $k \in X$.

Under \textbf{(A1)-(A2)}, it follows from the Lemma 3.3 of Hern\'{a}ndez-Lerma \cite{HL1} that for any
$(\mu,\nu) \in \SM_1 \times \SM_2$, the corresponding Markov chain
$\{X_t\}$ is uniformly ergodic with a unique invariant probability
measure $\eta[\mu,\nu] \in \PR(X)$, i.e.,
\begin{equation} \label{unierg}
||\hat{q}^t(\cdot|k,\mu,\nu) - \eta[\mu,\nu](\cdot)||_{\mbox{TV}} \leq 2 {\delta}^t,\ \forall k \in X,\ \forall t \in \NM,
\end{equation}
where $\hat{q}^t(\cdot|k,\mu,\nu)$ denotes the $t$-step transition kernel under $(\mu,\nu)$ starting from $k \in X$.

\begin{remark} \label{unifparambnd}
It is important to note here that the quantity $\delta$ on the r.h.s.\ of the above inequality (\ref{unierg}) is independent
of the strategy pair $(\mu,\nu)$ and \eqref{unierg} holds uniformly across all $(\mu,\nu) \in \SM_1 \times \SM_2$ (see also
condition C4 of Theorem 4(ii) in Federgruen et.\ al.\ \cite{FHT1}). Note that an assumption equivalent to \textbf{(A2)} was also
necessary for the corresponding zero-sum case as analyzed in Basu and Ghosh \cite{BaGh2}. Equivalent assumptions are quite common
in the corresponding literature for the one-controller case (see, e.g., Di Masi and Stettner \cite{DS1},\cite{DS1a} and \cite{DS2}, Hern\'{a}ndez-Hern\'{a}ndez and Marcus \cite{HM1}).
\end{remark}

%\begin{remark} \label{geomrecur}
\begin{definition} \label{geomrecur}
For any $B \subseteq X$, %and $(\mu,\nu) \in \SM_1 \times \SM_2$,
define the hitting time of the set $B$ by $\{X_t\}$ as
\begin{equation} \label{hittime}
\tau_B \defn \inf \{t\geq 0 : X_t \in B\}.
\end{equation}
%where $X^{\mu,\nu}_t$ is the process under $(\mu,\nu)$.
Denote $\tau_B = \tau_j$ if $B = \{j\}$. Similarly, define the first return time to $B$ by $\{X_t\}$ as
\begin{equation} \label{rettime}
\sigma_B \defn \inf\{t \geq 1 : X_t \in B\}.
\end{equation}
Denote $\sigma_B = \sigma_j$ if $B = \{j\}$.
\begin{comment}
Note that the Markov chain $\{X_t\}$ is said to be ``geometrically recurrent" under some $(\mu,\nu) \in \SM_1 \times \SM_2$ if there exists some $R_{\mu,\nu} > 1$ and a finite subset $C^{\mu,\nu}$ of $X$ such that $\sup_{k \in X}E^{\mu,\nu}_k[R_{\mu,\nu}^{\sigma_{C^{\mu,\nu}}}]  < \infty$.
\end{comment}
\end{definition}
%\end{remark}

\noi Now, Theorem 2.1 of Balaji and Meyn \cite{BaM1}, Theorem 16.0.2 of Meyn and Tweedie \cite{MT1} and Remark \ref{unifparambnd} together imply the following result:
\begin{proposition} \label{lpstblaprecur}
Assume \textbf{(A1)-(A2)}. Then the following three equivalent statements hold. %For $(\mu,\nu) \in \SM_1 \times \SM_2$,
\begin{enumerate}
\item[(a)] For every $A \in 2^{X}$, the chain $\{X_t\}$ is geometrically recurrent uniformly over all $(\mu,\nu) \in \SM_1 \times \SM_2$, i.e., there exists some $R_A > 1$ such that $\sup_{(\mu,\nu) \in \SM_1 \times \SM_2}\sup_{k \in X}E^{\mu,\nu}_k[R_A^{\sigma_A}] \leq B_A < \infty$. In particular, if $A = \{j\},\ j \in X$ we write the constants above as $R_j$ and $B_j$ respectively. Note that $B_A, B_j > 1$ by definition.
\item[(b)] For every $A \in 2^{X}$, there exists some $L_A >0$ such that $\sup_{(\mu,\nu) \in \SM_1 \times \SM_2}\sup_{k \in X}E^{\mu,\nu}_k[\sigma_A] \leq L_A < \infty$. In particular, if $A = \{j\},\ j \in X$ we write the constant above as $L_j$.
\item[(c)] There exist $\eta < 1, b < \infty$, a finite subset $C \subset X$ and a bounded function $V : X \ra [1,\infty)$ such that \[\sum_{j}V(j)q(j|i,\mu(i),\nu(i)) \leq \eta V(i) + b \mathbb{I}_C(i)\]
    where $\eta,b,C$ and $V$ are all independent of $(\mu,\nu) \in \SM_1 \times \SM_2$.
\end{enumerate}
\end{proposition}

\begin{remark} \label{examples}
Note that stochastic Lyapunov-type stability assumptions as in Proposition \ref{lpstblaprecur}(c) was also necessary for the
corresponding zero-sum stochastic differential game case as analyzed in Basu and Ghosh \cite{BaGh1}. See also Remark \ref{unifparambnd} for reference to the corresponding literature where equivalent assumptions are used. We also refer the reader
to Altman et.\ al.\ \cite{AltHorSpi1} for examples of such games on countable spaces with additive cost criteria where equivalent stability assumptions were made.
\end{remark}

Henceforth we fix $\theta \in (0,\Theta]$.
%Since $\theta$ is fixed, we do not show the explicit dependence on $\theta$ (where it is obvious) for notational convenience.
Given strategies $(\mu,\nu) \in \Pi_1 \times \Pi_2$, consider the following
dynamic programming (HJB) equations for the ergodic cost (\ref{riskergpayoff}) for players I and II (resp.):
\begin{equation} \label{HJIergupper}
e^{\theta\lambda_1 + V_1(\theta,k)} = \inf_{\xi \in \PR(U)} \left[\int_U \int_V e^{\theta r_1(k,u,v)}\left(\sum_{j \in X}e^{V_1(\theta,j)} q(j|k,u,v)\right) \xi(du)\nu_t[h_t](dv)\right],
\end{equation}
and
\begin{equation} \label{HJIerglower}
e^{\theta\lambda_2 + V_2(\theta,k)} = \inf_{\chi \in \PR(V)} \left[\int_U \int_V e^{\theta r_2(k,u,v)}\left(\sum_{j \in X}e^{V_2(\theta,j) }q(j|k,u,v)\right)\mu_t[h_t](du)\chi(dv)\right],
\end{equation}
where $\lambda_1,\lambda_2$ are scalars.

Under \textbf{(A1)}-\textbf{(A2)}, by Corollary 2.3 of Di Masi and Stettner \cite{DS1}, there exists at most one (up to an additive constant) bounded function $\hat{V}_1[\nu]$ and a unique constant $\hat{\lambda}_1[\nu]$ for which (\ref{HJIergupper}) is satisfied. If it does then, by Proposition 1.1 of Di Masi and Stettner \cite{DS1} (see also Theorem 2.1 of Hern\'{a}ndez-Hern\'{a}ndez and Marcus \cite{HM1}), we have
\begin{equation} \label{stochrep01erg}
\hat{\lambda}_1[\nu] = \inf_{\mu \in \Pi_1} \beta_1^{\mu,\nu}(\theta,k) = \inf_{\mu \in \SM_1}{\lim\sup}_{T \ra \infty}\frac{1}{\theta T} \ln E_k^{\mu,\nu}\left[e^{\theta \sum_{t=0}^{T-1} r_1(X_t,U_t,V_t)}\right]
\end{equation}
and similarly, there exists at most one bounded $\hat{V}_2[\mu]$ (modulo an additive constant) and a unique $\hat{\lambda}_2[\mu]$ satisfying (\ref{HJIerglower}) such that
\begin{equation} \label{stochrep02erg}
\hat{\lambda}_2[\mu] = \inf_{\nu \in \Pi_2} \beta_2^{\mu,\nu}(\theta,k) = \inf_{\nu \in \SM_2}{\lim\sup}_{T \ra \infty}\frac{1}{\theta T} \ln E_k^{\mu,\nu}\left[e^{\theta \sum_{t=0}^{T-1} r_2(X_t,U_t,V_t)}\right].
\end{equation}
%with $\{X^{\mu,\nu}_t,U^{\mu,\nu}_t,V^{\mu,\nu}_t : t \in \NM_0\}$ being the corresponding processes generated by $(\mu,\nu) \in \Pi_1 \times \Pi_2$. %Moreover, by the continuity of $q$ and the compactness of $\PR(U)$ and $\PR(V)$,
Hence there exist measurable (minimizing) selectors (see Bene\v{s} \cite{Be1})
\[(\hat{\mu}[\nu],\hat{\nu}[\mu]) : X \mapsto \PR(U) \times \PR(V)\] satisfying
\begin{eqnarray} \label{discminsel1erg}
\inf_{\xi \in \PR(U)} \left[\int_U \int_V e^{\theta r_1(k,u,v)}\left(\sum_{j \in X}e^{\hat{V}_1[\nu](\theta,j)} q(j|k,u,v)\right)\xi(du)\nu_t[h_t](dv)\right] \nonumber \\
= \int_U \int_V e^{\theta r_1(k,u,v)} \left(\sum_{j \in X} e^{\hat{V}_1[\nu](\theta,j)}q(j|k,u,v)\right)\hat{\mu}[\nu]\left[du|k\right]\nu_t[h_t](dv),
\end{eqnarray}
and
\begin{eqnarray} \label{discminsel2erg}
\inf_{\chi \in \PR(V)} \left[\int_U \int_V e^{\theta r_2(k,u,v)} \left(\sum_{j \in X} e^{\hat{V}_2[\mu](\theta,j)} q(j|k,u,v)\right)\mu_t[h_t](du)\chi(dv)\right] \nonumber \\
= \int_U \int_V e^{\theta r_2(k,u,v)} \left(\sum_{j \in X} e^{\hat{V}_2[\mu](\theta,j)}q(j|k,u,v)\right)\mu_t[h_t](du)\hat{\nu}[\mu]\left[dv|k\right].
\end{eqnarray}

Hence, by Proposition 1.1 of Di Masi and Stettner \cite{DS1} (see also Theorem 2.1 of Hern\'{a}ndez-Hern\'{a}ndez and Marcus \cite{HM1}), given $(\mu,\nu) \in \Pi_1 \times \Pi_2$, if there exist bounded solutions to (\ref{HJIergupper}) and (\ref{HJIerglower}) then the minimizing
%strategies $\{u_t^*\} \in \SM_1,\ \{v_t^*\} \in \SM_2$ in (\ref{stochrep01erg}) and (\ref{stochrep02erg}) (resp.) are given by (see Proposition 1.1 of Di Masi and Stettner \cite{DS1})
selector $\hat{\mu}[\nu]$ and $\hat{\nu}[\mu]$ in (\ref{discminsel1erg}) and (\ref{discminsel2erg}) resp.
%\begin{equation} \label{minseq01erg}
%u^*_t = \hat{\mu}[\nu](X_t),
%\end{equation}
%and
%\begin{equation} \label{minseq02erg}
%v^*_t = \hat{\nu}[\mu](X_t),
%\end{equation}
 generates stationary (see \eqref{stochrep01erg} and \eqref{stochrep02erg}) optimal response strategy of player I (resp.\ player II) when the other's strategy fixed. Hence, without loss of generality, we can effectively consider the strategy pairs $(\mu,\nu) \in \SM_1 \times \SM_2$ for further analysis of this game. Thus if we consider the minimizing selectors $\hat{\mu}[\nu]$ and $\hat{\nu}[\mu]$ in (\ref{discminsel1erg}) and (\ref{discminsel2erg}) respectively, we can define the following point-to-set maps $\HL_i : \SM_j \mapsto 2^{\SM_i},\ i,j=1,2,\ i\neq j$ as follows:
\begin{eqnarray} \label{optresponseerg01}
\HL_1(\nu) \defn \{\{\mu^*[\nu]\} : \mu^*[\nu]\ \mbox{satisfies}\ (\ref{discminsel1erg})\},
\end{eqnarray}
and
\begin{eqnarray} \label{optresponseerg02}
\HL_2(\mu) \defn \{\{\nu^*[\mu]\} : \nu^*[\mu]\ \mbox{satisfies}\ (\ref{discminsel2erg})\}.
\end{eqnarray}

Define the map $\HL \equiv \HL_1 \otimes \HL_2 : \SM_1 \times \SM_2 \mapsto 2^{\SM_1 \times \SM_2}$ as
\begin{eqnarray} \label{optresponseerg}
\HL(\mu,\nu) \defn \{(\{\mu^*[\nu]\},\{\nu^*[\mu]\}) : \mu^*[\nu]\ \mbox{satisfies}\ (\ref{discminsel1erg})\ \mbox{and}\ \nu^*[\mu]\ \mbox{satisfies}\ (\ref{discminsel2erg})\}.
\end{eqnarray}

A priori, it is not clear that the sets $\HL_1, \HL_2$ and hence $\HL$ are nonempty. We now proceed to prove the existence of unique bounded solutions to (\ref{HJIergupper}) and (\ref{HJIerglower}). Once this result is established, the nonemptyness of the above sets follow automatically. Then, using arguments similar to Theorem \ref{existexpvalue}, we prove the existence of Nash equilibria by showing that the map $\HL$ has a fixed point $(\tilde{\mu},\tilde{\nu})$. To this end, we make the following definitions (see Balaji and Meyn \cite{BaM1} and Borkar and Meyn \cite{BM1}).% and prove the following supporting lemmata.

\begin{definition} \label{balgurudef}
Given an arbitrarily fixed state $0 \in X$ and any strategy $(\mu,\nu) \in \SM_1 \times \SM_2$, %with $\{X^{\mu,\nu}_t,U^{\mu,\nu}_t,V^{\mu,\nu}_t : t \in \NM_0\}$ being the corresponding processes generated by $(\mu,\nu)$,
let
\begin{equation} \label{generalgpe01}
\Lambda_1[\nu](\mu) \defn \inf\left\{\Lambda \in \R : E_0^{\mu,\nu}\left[e^{\theta\sum_{t=0}^{\sigma_0 - 1}(r_1(X_t,U_t,V_t) - \Lambda)}\right] \leq 1\right\},
\end{equation}
and
\begin{equation} \label{generalgpe02}
\Lambda_2[\mu](\nu) \defn \inf\left\{\Lambda \in \R : E_0^{\mu,\nu}\left[e^{\theta\sum_{t=0}^{\sigma_0 - 1}(r_2(X_t,U_t,V_t) - \Lambda)}\right] \leq 1\right\}.
\end{equation}

Also define
\begin{eqnarray}
\Lambda_1^*[\nu] &\defn& \inf_{\mu \in \SM_1} \Lambda_1[\nu](\mu), \label{mingpe01} \\
h_1^*[\nu](k) &\defn& \inf_{\mu \in \SM_1}E_k^{\mu,\nu}\left[e^{\theta\sum_{t=0}^{\tau_0}(r_1(X_t,U_t,V_t) - \Lambda_1^*[\nu])}\right], \label{relval01} \\
w_1^*[\nu](k) &\defn& {\arg\min}_{\xi \in \PR(U)} \left[\int_U \int_V e^{\theta r_1(k,u,v)}\left(\sum_{j\in X} h_1^*[\nu](j)q(j|k,u,v)\right)\xi(du)\nu\left[dv|k\right]\right], \label{mpearg01}
\end{eqnarray}
where $w_1^*[\nu]$ is any minimizer in (\ref{mpearg01}). Similarly, let
\begin{eqnarray}
\Lambda_2^*[\mu] &\defn& \inf_{\nu \in \SM_2} \Lambda_2[\mu](\nu), \label{mingpe02} \\
h_2^*[\mu](k) &\defn& \inf_{\nu \in \SM_2}E_k^{\mu,\nu}\left[e^{\theta\sum_{t=0}^{\tau_0}(r_2(X_t,U_t,V_t) - \Lambda_2^*[\mu])}\right], \label{relval02} \\
w_2^*[\mu](k) &\defn& {\arg\min}_{\chi \in \PR(V)} \left[\int_U \int_V e^{\theta r_2(k,u,v)}\left(\sum_{j\in X} h_2^*[\mu](j)q(j|k,u,v)\right)\mu\left[du|k\right]\chi(dv)\right] \label{mpearg02}
\end{eqnarray}
where $w_2^*[\mu]$ is any minimizer in (\ref{mpearg02}).
\end{definition}

Note that since $\theta$ is fixed we have suppressed the explicit dependence on $\theta$ in these definitions for notational convenience. A priori, it is not obvious that the quantities defined above are finite. The following Proposition \ref{eqmingpeergcost} and Lemmata \ref{finrelvalmpeineq} and \ref{relvalfnbndlem} settle this issue.

First, we %choose one special fixed but arbitrary state $0 \in X$ and
make a ``small cost" assumption which shall hold for the rest of this paper.  \\ %is needed for proving the uniqueness of the solutions to (\ref{HJIergupper}) and (\ref{HJIerglower}):\\
\noi \textbf{(A3)} $||r_i||_\infty \leq \frac{\ln R_0}{3 \Theta},\ i =1,2$,\\
where $R_0 > 1$ is as in part (a) of Proposition \ref{lpstblaprecur}.

\begin{remark} \label{smallcost}
Note that a similar assumption was also
necessary for the corresponding zero-sum stochastic differential game case as analyzed in Basu and Ghosh \cite{BaGh1}.
\end{remark}

\begin{comment} %%%%%%%%%%%%%% SERFOZO %%%%%%%%%%%%%%%
\begin{definition} \label{integrable}
Given $\varsigma, \varsigma_n\ (n \geq 1)$ nonnegative $\sigma$-finite measures on a measurable space $(\Omega,\F)$ and $f,f_n\ (n \geq 1)$
measurable functions from $(\Omega,\F)$ to $(\R, \B(\R))$, we call the sequence $\{f_n\}$ ``uniformly $\{\varsigma_n\}$-integrable" if
\begin{equation} \label{ui}
\lim_{z \ra \infty} \sup_n \int_{|f_n| > z} |f_n| d\varsigma_n = 0,
\end{equation}
and ``tightly $\{\varsigma_n\}$-integrable" if
\begin{equation} \label{ti}
\inf_{\{B_n\} \in \mathfrak{S}} \sup_n \int_{B_n^c} |f_n| d\varsigma_n = 0
\end{equation}
where $\mathfrak{S} \defn \{\{B_n\}:B_n \in \F\ \mbox{and}\ \sup_n \ \varsigma_n(B_n) < \infty\}$. We say that $f_n \ra f$ in $\varsigma_n$-measure if, for any $\ep > 0$,
\begin{equation} \label{measconv}
\varsigma_n\{\omega: |f_n(\omega)-f(\omega)| \geq \ep\} \stackrel{n \ra \infty}{\lra} 0.
\end{equation}
\end{definition}

\begin{remark} \label{equivti}
Note that (\ref{ti}) holds when $\sup_n \varsigma_n(\Omega) < \infty$, in particular, when $\varsigma_n$ are probability measures (see p.\ 383 of Serfozo \cite{Ser} for details).
\end{remark}
\end{comment} %%%%%%%%%%%%%%%%%%% SERFOZO %%%%%%%%%%%%%%%%%%%%%

Now we state and prove the following important proposition.

\begin{proposition} \label{eqmingpeergcost}
Under \textbf{(A1)}-\textbf{(A3)} and for any $(\mu,\nu) \in \SM_1 \times \SM_2$,
\begin{equation} \label{eqmingpeergcostbnd01}
||r_1||_{\infty} \geq  |\Lambda_1^*[\nu]|,
\end{equation}
and
\begin{equation} \label{eqmingpeergcostbnd02}
||r_2||_{\infty} \geq |\Lambda_2^*[\mu]|.
\end{equation}
Moreover, $\Lambda_1^*[\nu]$ (resp.\ $\Lambda_2^*[\mu]$) is continuous in $\nu$ (resp.\ $\mu$).
\end{proposition}

\noi \textbf{Proof:} We prove it for player I. The result for player II follows analogously. From (\ref{generalgpe01}) we obtain, by Jensen's inequality
\begin{eqnarray*}
1 &&\geq E_0^{\mu,\nu}\left[e^{\theta\sum_{t=0}^{\sigma_0 - 1}(r_1(X_t,U_t,V_t) - \Lambda_1[\nu](\mu))}\right] \geq e^{-\theta \Lambda_1[\nu](\mu)E_0^{\mu,\nu}[\sigma_0]} e^{\theta E_0^{\mu,\nu}\left[\sum_{t=0}^{\sigma_0 - 1}(r_1(X_t,U_t,V_t)\right]} \\ && \geq e^{-\theta \Lambda_1[\nu](\mu)E_0^{\mu,\nu}[\sigma_0]} e^{-\theta||r_1||_{\infty} E_0^{\mu,\nu}[\sigma_0]}
\end{eqnarray*}
which implies
\begin{eqnarray*}
e^{\theta||r_1||_{\infty} E_0^{\mu,\nu}[\sigma_0]} \geq e^{-\theta \Lambda_1[\nu](\mu)E_0^{\mu,\nu}[\sigma_0]}.
\end{eqnarray*}
Under \textbf{(A1)}-\textbf{(A2)}, $0 < E_0^{\mu,\nu}[\sigma_0] \leq L_0 < \infty$ (see part (b) of Proposition \ref{lpstblaprecur}) and hence, by taking $\log$ on both sides of the last inequality, we get \[-||r_1||_{\infty} \leq  \Lambda_1[\nu](\mu).\]
Now, by (\ref{generalgpe01}), for any $\ep > 0$, we have
\begin{eqnarray*}
E_0^{\mu,\nu}\left[e^{\theta\sum_{t=0}^{\sigma_0 - 1}(r_1(X_t,U_t,V_t) - \Lambda_1[\nu](\mu) + \ep)}\right] > 1
\end{eqnarray*}
implying
\begin{eqnarray*}
E_0^{\mu,\nu}\left[e^{\theta\sum_{t=0}^{\sigma_0 - 1}(||r_1||_{\infty} - \Lambda_1[\nu](\mu) + \ep)}\right] > 1
\end{eqnarray*}
which, in turn, implies
\begin{eqnarray*}
E_0^{\mu,\nu}\left[e^{\theta(||r_1||_{\infty} - \Lambda_1[\nu](\mu) + \ep)\sigma_0}\right] > 1.
\end{eqnarray*}
Since this is true for any $\ep > 0$, we have in the limit $\ep \downarrow 0$,
\begin{eqnarray*}
E_0^{\mu,\nu}\left[e^{\theta(||r_1||_{\infty} - \Lambda_1[\nu](\mu))\sigma_0}\right] \geq 1
\end{eqnarray*}
implying
\begin{eqnarray*}
||r_1||_{\infty} - \Lambda_1[\nu](\mu) \geq 0.
\end{eqnarray*}
Hence we have
\begin{equation} \label{lambdabnd}
-||r_1||_{\infty} \leq \Lambda_1[\nu](\mu) \leq ||r_1||_{\infty}.
\end{equation}
The inequality (\ref{eqmingpeergcostbnd01}) follows by taking the $\inf$ over $\mu \in \SM_1$. We now prove the joint continuity in $(\mu,\nu)$ of $\Lambda_1[\nu](\mu)$. Consider a sequence $\{(\mu^{(n)},\nu^{(n)})\} \stackrel{n \ra \infty}{\lra} (\mu^{(\infty)},\nu^{(\infty)})$ in $\SM_1 \times \SM_2$. Let $\{X_t,t\geq 0\}$ denote the chain starting at $0$. By Proposition \ref{contlawstrat}, \[\LA^{\mu^{(n)},\nu^{(n)}}\{X_t, t\geq 0\} \stackrel{n \ra \infty}{\lra} \LA^{\mu^{(\infty)},\nu^{(\infty)}}\{X_t, t \geq 0\}.\] Then, by Skorohod's Theorem (see Theorem 2.2.2 of Borkar \cite{B1}),
there exists some augmentation $(\tilde{H}_{\infty},\B(\tilde{H}_{\infty}),\tilde{\PR})$ of the canonical space $(H_{\infty},\B(H_{\infty}))$ on which \[\{X^{\mu^{(n)},\nu^{(n)}}_t, t\geq 0\} \stackrel{\tilde{P}\mbox{-}a.s}{\lra} \{X^{\mu^{(\infty)},\nu^{(\infty)}}_t, t\geq 0\}\] where
$\{X^{\mu^{(n)},\nu^{(n)}}_t, t\geq 0\}, n = 1,2,\ldots,\infty$ denotes the Markov chain under $(\mu^{(n)},\nu^{(n)})$ on this augmented space starting at $0$. Following (\ref{rettime}), let $\sigma_0^{(n)} \defn \inf\{t \geq 1 : X^{\mu^{(n)},\nu^{(n)}}_t = 0\},\ n=1,2,\ldots,\infty$. Then, obviously $\tilde{P}$-a.s.\[\sigma_0^{(n)} \stackrel{n \ra \infty}{\lra} \sigma_0^{(\infty)}.\] Hence we get the joint continuity under $(\mu,\nu)$ of $\sigma_0$.
\begin{comment}
This means that for all $\omega \in \tilde{H}_{\infty} \setminus \A$ with $\tilde{\PR}(\A) = 0$, there exists some finite $M(\omega) \geq 1$ such that for all $n \geq M(\omega)$ and for all $t \geq 0$, $X^{\mu^{(n)},\nu^{(n)}}_t (\omega) = X^{\mu^{(\infty)},\nu^{(\infty)}}_t (\omega)$ as the state space is discrete.  Following (\ref{rettime}), let $\sigma_B^{(n)} \defn \inf\{t \geq 1 : X^{\mu^{(n)},\nu^{(n)}}_t \in B\},\ n=1,2,\ldots,\infty$ for any $B \subseteq X$. Now suppose for some $A \subset X,\ \sigma_A^{(n)} \stackrel{a.s.}{\nrightarrow} \sigma_A^{(\infty)}$ in the discrete topology. Then there exist some $\mathcal{N} \in \B(\tilde{H}_{\infty})$ with $\tilde{\PR}(\mathcal{N}) > 0$ such that for each $\omega \in \mathcal{N}$ there exists infinitely many $n_k (\omega),k=1,2,\ldots$ with $\inf\{t \geq 1 : X^{\mu^{(n_k (\omega))},\nu^{(n_k (\omega))}}_t (\omega) \in A\} \neq \inf\{t \geq 1 : X^{\mu^{(\infty)},\nu^{(\infty)}}_t (\omega) \in A\}$, i.e., there exists some finite $t(\omega) \geq 1$ with $X^{\mu^{(n_k (\omega))},\nu^{(n_k (\omega))}}_{t(\omega)} (\omega) \neq X^{\mu^{(\infty)},\nu^{(\infty)}}_{t(\omega)} (\omega)$. Note that the finiteness of $t(\omega)$ follows from the fact that $E_0^{\mu,\nu}[\sigma_0^{\mu,\nu}] \leq L_0$ for any $(\mu,\nu) \in \SM_1 \times \SM_2$ under \textbf{(A1)}-\textbf{(A2)} (see part (b) of Proposition \ref{lpstblaprecur}). Thus we get a contradiction implying the joint continuity under $(\mu,\nu)$ of $\sigma^{\mu,\nu}_A$, in particular, $\sigma^{\mu,\nu}_0$.
\end{comment}
Hence, for any $0 < \ep \leq ||r_1||_{\infty}$ and along a subsequence $\{n_k\} \subset \{n\}$ with $\{(\mu^{(n_k)},\nu^{(n_k)})\} \stackrel{n_k \ra \infty}{\lra} (\mu^{(\infty)},\nu^{(\infty)})$, we have by (\ref{lambdabnd}) and the joint continuity of $\sigma_0$ under $(\mu,\nu)$
\begin{eqnarray} \label{lawgpelhs}
f_{n_k} \equiv e^{\theta\sum_{t=0}^{\sigma_0^{(n_k)} - 1}\left(r_1(X_t^{\mu^{(n_k)},\nu^{(n_k)}},U_t^{\mu^{(n_k)},\nu^{(n_k)}},V_t^{\mu^{(n_k)},\nu^{(n_k)}}) - \Lambda_1[\nu^{(n_k)}](\mu^{(n_k)}) + \ep\right)} \nonumber \\ \stackrel{\tilde{P}\mbox{-}a.s.}{\lra}
e^{\theta\sum_{t=0}^{\sigma_0^{(\infty)} - 1}\left(r_1(X_t^{\mu^{(\infty)},\nu^{(\infty)}},U_t^{\mu^{(\infty)},\nu^{(\infty)}},V_t^{\mu^{(\infty)},\nu^{(\infty)}}) - \Lambda + \ep\right)} \equiv f_{\infty},
\end{eqnarray}
where $|\Lambda| \leq ||r_1||_{\infty}$ and $U_t^{\mu^{(n)},\nu^{(n)}},V_t^{\mu^{(n)},\nu^{(n)}}$ are the corresponding control sequences under $(\mu^{(n)},\nu^{(n)}),\ n=1,2,\ldots,\infty$ on this augmented space. Then
\begin{equation} \label{lawgpelhsbnd}
0 < f_{n_k} = |f_{n_k}| \leq g_{n_k} \defn e^{\theta(2||r_1||_{\infty} + \ep)\sigma_0^{(n_k)}} \stackrel{\tilde{P}\mbox{-}a.s.}{\lra} e^{\theta(2||r_1||_{\infty} + \ep)\sigma_0^{(\infty)}} \defn g_{\infty} \geq |f_{\infty}| = f_{\infty}.
\end{equation}
%Note that by Remark \ref{equivti}, $\{g_{n_k}\}$ is tightly $\left\{P^{\mu^{(n_k)},\nu^{(n_k)}}_0\right\}$-integrable.
Moreover, by \textbf{(A3)} and part (a) of Proposition \ref{lpstblaprecur}, we have
\begin{equation} \label{uiequivcond}
\sup_{n_k}\tilde{E}\left[g_{n_k}\right] \leq \sup_{n_k}\tilde{E}\left[e^{3\Theta||r_1||_{\infty}\sigma_0^{(n_k)}}\right] \leq \sup_{n_k}\tilde{E}\left[R_0^{\sigma_0^{(n_k)}}\right] \leq B_0 < \infty,
\end{equation}
which, in particular, implies that $\{\sigma_0^{(n_k)}\}$ are uniformly integrable.
Given any $\varepsilon > 0$ it follows from (\ref{uiequivcond}) that there exists some large enough constant $K \equiv K(\varepsilon,B_0) > 0$ such that
\begin{equation} \label{uistep}
\sup_{n_k}\tilde{E}\left[g_{n_k} \mathbb{I}_{\{g_{n_k} \geq K\}}\right] \leq \varepsilon
\end{equation}
implying $\{g_{n_k}\}$ are uniformly integrable.

\begin{comment}
Now let there be any $\{A_{n_k}\} \subset \B(\tilde{H}_{\infty})$ with $\sup_{n_k} P^{\mu^{(n_k)},\nu^{(n_k)}}_0\left(A_{n_k}\right) < \frac{\varepsilon}{2 K}$. Then, from (\ref{uistep}), we have
\begin{eqnarray} \label{uifinal}
\sup_{n_k}E_0^{\mu^{(n_k)},\nu^{(n_k)}}\left[g_{n_k}\mathbb{I}_{A_{n_k}}\right] && \leq \sup_{n_k}E_0^{\mu^{(n_k)},\nu^{(n_k)}}\left[g_{n_k}\mathbb{I}_{A_{n_k} \bigcap \{g_{n_k} \leq K\}}\right] + \sup_{n_k}E_0^{\mu^{(n_k)},\nu^{(n_k)}}\left[g_{n_k}\mathbb{I}_{A_{n_k}\bigcap \{g_{n_k} \geq K\}}\right] \nonumber \\
&&\leq K \sup_{n_k} P^{\mu^{(n_k)},\nu^{(n_k)}}_0\left(A_{n_k}\right) + \sup_{n_k}E_0^{\mu^{(n_k)},\nu^{(n_k)}}\left[g_{n_k}\mathbb{I}_{\{g_{n_k} \geq K\}}\right] \nonumber \\ &&\leq K \frac{\varepsilon}{2K} + \frac{\varepsilon}{2} \leq \varepsilon.
\end{eqnarray}
\end{comment}

%Hence, by Lemma 2.5 of Serfozo \cite{Ser}, $\{g_{n_k}\}$ is uniformly and tightly $\left\{P^{\mu^{(n_k)},\nu^{(n_k)}}_0\right\}$-integrable.

It follows by Theorem 6.18(iii) of Klenke \cite{Kl1} %Remark 2.6(b) of Serfozo \cite{Ser}
and (\ref{lawgpelhsbnd}) that $\{f_{n_k}\}$ is uniformly integrable.

%Moreover, (\ref{lawgpelhs}) implies $f_{n_k} \stackrel{n \ra \infty}{\lra} f_{\infty}$ in  $P^{\mu^{(n_k)},\nu^{(n_k)}}_0$-measure, i.e., (\ref{measconv}) holds.

Also, by \textbf{(A3)} and part (a) of Proposition \ref{lpstblaprecur}, we have
%Note that, by the joint continuity under $(\mu,\nu)$ of $\sigma_0$, Proposition \ref{contlawstrat} and \textbf{(A3)}, we have
\begin{equation} \label{finitelimstrat}
\tilde{E}\left[g_{\infty}\right] \leq \tilde{E}\left[e^{3\Theta||r_1||_{\infty}\sigma_0^{(\infty)}}\right] \leq \tilde{E}\left[R_0^{\sigma_0^{(\infty)}}\right] \leq B_0 < \infty.
\end{equation}
Hence, by Corollary 1.3.1 of Borkar \cite{B1}, we get %Theorem 2.8 of Serfozo \cite{Ser}, we get %by Proposition 18, Chapter 11 of Royden \cite{Ry}, we get
\begin{eqnarray} \label{lawgpelhsfin}
&&\lim_{n_k \ra \infty}\tilde{E}\left[e^{\theta\sum_{t=0}^{\sigma_0^{(n_k)} - 1}\left(r_1(X_t^{\mu^{(n_k)},\nu^{(n_k)}},U_t^{\mu^{(n_k)},\nu^{(n_k)}},V_t^{\mu^{(n_k)},\nu^{(n_k)}}) - \Lambda_1[\nu^{(n_k)}](\mu^{(n_k)}) + \ep\right)}\right] \nonumber \\ &&= \tilde{E}\left[e^{\theta\sum_{t=0}^{\sigma_0^{(\infty)} - 1}\left(r_1(X_t^{\mu^{(\infty)},\nu^{(\infty)}},U_t^{\mu^{(\infty)},\nu^{(\infty)}},V_t^{\mu^{(\infty)},\nu^{(\infty)}}) - \Lambda + \ep\right)}\right].
\end{eqnarray}
Similarly we can show that
\begin{eqnarray} \label{lawgperhsfin}
&&\lim_{n_k \ra \infty}\tilde{E}\left[e^{\theta\sum_{t=0}^{\sigma_0^{(n_k)} - 1}\left(r_1(X_t^{\mu^{(n_k)},\nu^{(n_k)}},U_t^{\mu^{(n_k)},\nu^{(n_k)}},V_t^{\mu^{(n_k)},\nu^{(n_k)}}) - \Lambda_1[\nu^{(n_k)}](\mu^{(n_k)})\right)}\right] \nonumber \\ &&= \tilde{E}\left[e^{\theta\sum_{t=0}^{\sigma_0^{(\infty)} - 1}\left(r_1(X_t^{\mu^{(\infty)},\nu^{(\infty)}},U_t^{\mu^{(\infty)},\nu^{(\infty)}},V_t^{\mu^{(\infty)},\nu^{(\infty)}}) - \Lambda\right)}\right].
\end{eqnarray}
Now from (\ref{generalgpe01}) we get that for any $0 < \ep \leq ||r_1||_{\infty}$,
\begin{eqnarray} \label{jntcntlawgpe}
&&\tilde{E}\left[e^{\theta\sum_{t=0}^{\sigma_0^{(n)} - 1}\left(r_1(X_t^{\mu^{(n)},\nu^{(n)}},U_t^{\mu^{(n)},\nu^{(n)}},V_t^{\mu^{(n)},\nu^{(n)}}) - \Lambda_1[\nu^{(n)}](\mu^{(n)}) + \ep\right)}\right] \geq 1 \nonumber \\
&& \geq \tilde{E}\left[e^{\theta\sum_{t=0}^{\sigma_0^{(n)} - 1}\left(r_1(X_t^{\mu^{(n)},\nu^{(n)}},U_t^{\mu^{(n)},\nu^{(n)}},V_t^{\mu^{(n)},\nu^{(n)}}) - \Lambda_1[\nu^{(n)}](\mu^{(n)})\right)}\right].
\end{eqnarray}
Hence, taking $\lim_{n_k \ra \infty}$ on both sides of (\ref{jntcntlawgpe}), we have by \eqref{lawgperhsfin} and \eqref{jntcntlawgpe}
\begin{eqnarray} \label{jntcntlawgpeinf}
&&\tilde{E}\left[e^{\theta\sum_{t=0}^{\sigma_0^{(\infty)} - 1}\left(r_1(X_t^{\mu^{(\infty)},\nu^{(\infty)}},U_t^{\mu^{(\infty)},\nu^{(\infty)}},V_t^{\mu^{(\infty)},\nu^{(\infty)}}) - \Lambda + \ep\right)}\right] \geq 1 \nonumber \\ &&\geq \tilde{E}\left[e^{\theta\sum_{t=0}^{\sigma_0^{(\infty)} - 1}\left(r_1(X_t^{\mu^{(\infty)},\nu^{(\infty)}},U_t^{\mu^{(\infty)},\nu^{(\infty)}},V_t^{\mu^{(\infty)},\nu^{(\infty)}}) - \Lambda\right)}\right],
\end{eqnarray}
for some $\Lambda$ satisfying (\ref{lambdabnd}) and all $0 < \ep \leq ||r_1||_{\infty}$. This implies, by (\ref{generalgpe01}), that $\Lambda = \Lambda_1[\nu^{(\infty)}](\mu^{(\infty)})$. Hence $\Lambda_1[\nu](\mu)$ is jointly continuous under $(\mu,\nu)$. The continuity of $\Lambda_1^*[\nu]$ with respect to $\nu$ now follows from (\ref{lambdabnd}) since $\SM_1$ is compact. \hfill $\Box$

\begin{remark} \label{remriskfacbnd}
Note that by (\ref{lambdabnd}) above the quantities $\Lambda_1[\nu](\mu)$ and $\Lambda_2[\mu](\nu)$ defined in (\ref{generalgpe01}) and (\ref{generalgpe02})
respectively are uniformly bounded for all $\theta \in (0,\Theta]$. This fact shall be implicitly used in the proof of Theorem \ref{relvalfnrepmpeuniq}.
\end{remark}

Now we state and prove two useful lemmata.

\begin{lemma} \label{finrelvalmpeineq}
Under \textbf{(A1)}-\textbf{(A2)} and for any $(\mu,\nu) \in \SM_1 \times \SM_2$,
\begin{enumerate}
\item The multiplicative Poisson inequalities holds for all $k \in X$, i.e.,
\begin{eqnarray} \label{mpineq}
&&\int_U \int_V e^{\theta r_1(k,u,v)}\left(\sum_{j \in X}h_1^*[\nu](j) q(j|k,u,v) \right) w_1^*[\nu](k)(du)\nu\left[dv|k\right] \leq e^{\theta\Lambda_1^*[\nu]} h_1^*[\nu](k), \\
&& \int_U \int_V e^{\theta r_2(k,u,v)}\left(\sum_{j \in X}h_2^*[\mu](j) q(j|k,u,v)\right)\mu\left[du|k\right]w_2^*[\mu](k)(dv) \leq e^{\theta\Lambda_2^*[\mu]} h_2^*[\mu](k).
\end{eqnarray}
\item $h_1^*[\nu](k) < \infty,\ h_2^*[\mu](k) < \infty,\ \forall k \in X$.
\end{enumerate}
\end{lemma}

\noi \textbf{Proof:} We prove it for player I when player II uses $\nu \in \SM_2$. The other case is analogous. %From Lemma \ref{eqmingpeergcost} it follows that $\Lambda_1^*[\nu](\theta) < \inft$ and there exists a $\SM_1 \ni \tilde{\mu}[\nu] \equiv (\tilde{w}_1[\nu],\tilde{w}_1[\nu],\ldots,) \in \left((\PR(U))^X\right)^{\mathbb{Z}_+}$ such that $\Lambda_1^*[\nu](\theta) = \Lambda_1[\nu](\theta,\tilde{\mu}[\nu])$.
For each $k \neq 0$,
\begin{eqnarray*}
&& e^{-\theta\Lambda_1^*[\nu]} \int_V \int_U e^{\theta r_1(k,u,v)} \left(\sum_{j \in X} h_1^*[\nu](j) q(j|k,u,v)\right)w_1^*[\nu](k)(du)\nu\left[dv|k\right] = \\ && \inf_{\xi \in \PR(U)} \int_V \int_U e^{\theta (r_1(k,u,v) - \Lambda_1^*[\nu])} \left(\sum_{j \in X} q(j|k,u,v) \left\{\inf_{\mu \in \SM_1}E_j^{\mu,\nu}\left[e^{\theta\sum_{t=0}^{\tau_0-1} \left(r_1(X_t,U_t,V_t) - \Lambda_1^*[\nu]\right)}\right]\right\}\right) \xi(du)\nu\left[dv|k\right]  \\ && = h_1^*[\nu](k)\ \mbox{(by one-step backward recursion)},
\end{eqnarray*}
while for $k = 0$, we have
\begin{eqnarray}
\int_V \int_U \left(\sum_{j \in X} h_1^*[\nu](j) q(j|0,u,v)\right) w_1^*[\nu](0)(du)\nu\left[dv|0\right] = %\nonumber \\ &&
\inf_{\mu \in \SM_1}E_0^{\mu,\nu}\left[e^{\theta\sum_{t=0}^{\sigma_0-1} (r_1(X_t,U_t,V_t) - \Lambda_1^*[\nu])}\right] \leq 1, \nonumber \\ \label{geomrecurcond} %\nonumber \\ && = E_0^{\tilde{\mu}[\nu],\nu}\left[e^{\theta\sum_{t=0}^{\sigma_0-1} (r_1(X_t) - \Lambda_1[\nu](\theta,\tilde{\mu}[\nu]))}\right] \leq 1, \label{geomrecurcond}
\end{eqnarray}
where %the last equality follows from Proposition \ref{contlawstrat} and
the inequality follows from (\ref{generalgpe01}). It follows that (\ref{mpineq}) holds. From (\ref{mpineq}) it follows exactly as in the proof of Lemma 15.2.2(i) of Meyn and Tweedie \cite{MT1} that the set $S \defn \{k \in X : h_1^*[\nu](k) < \infty\}$ absorbing, i.e., $\sum_{j \in S} \int_V \int_U  q(j|k,u,v) w_1^*[\nu](k)(du)\nu\left[dv|k\right] = 1$ for all $k \in S$. Since it is non-empty as $0 \in S$ and the chain $\{X_t\}$ is irreducible under the stationary strategy $(w_1^*[\nu],\nu)$  %($\{w_1^*[\nu]\}$ is the stationary strategy generated by $w_1^*[\nu]$)
by assumption \textbf{(A1)}, we get by Proposition 4.2.3 of Meyn and Tweedie \cite{MT1} that it is full, i.e. $S = X$ establishing part $2$. \hfill $\Box$

Note that, under the additional assumption \textbf{(A3)}, we can actually show a tighter result than part 2 of Lemma \ref{finrelvalmpeineq} as given below.
%Now we prove state and prove two important lemmata.

\begin{lemma} \label{relvalfnbndlem}
Under \textbf{(A1)} - \textbf{(A3)} and for any $(\mu,\nu) \in \SM_1 \times \SM_2$,
\begin{equation} \label{relvalfnbnd}
\frac{1}{R_0 B_0} \leq h_1^*[\nu](k) \leq R_0 B_0,\ \frac{1}{R_0 B_0} \leq \ h_2^*[\mu](k) \leq R_0 B_0,\ \forall k \in X.
\end{equation}
Moreover, $h_1^*[\nu](k)$ (resp.\ $h_2^*[\mu](k)$) is continuous in $\nu$ (resp.\ $\mu$) for all $k \in X$.
\end{lemma}

\noi \textbf{Proof:} We prove it for $h_1^*[\nu]$. The proof is exactly similar for $h_2^*[\mu]$. Under the assumptions of this lemma, we have by (\ref{relval01}) and part (a) of Proposition \ref{lpstblaprecur} %\ref{eqmingpeergcost},
\begin{eqnarray} \label{relvalbnd00}
0 &&< E_k^{\mu,\nu}\left[e^{\theta\sum_{t=0}^{\tau_0}(r_1(X_t,U_t,V_t) - \Lambda_1^*[\nu])}\right] =  E_k^{\mu,\nu}\left[e^{\theta\sum_{t=0}^{\sigma_0}(r_1(X_t,U_t,V_t) - \Lambda_1^*[\nu])}\right] \nonumber \\ &&\leq E_k^{\mu,\nu}\left[e^{2\theta||r_1||_\infty (\sigma_0+1)}\right] %\nonumber\\ &&
\leq E_k^{\mu,\nu}\left[e^{3\Theta||r_1||_\infty (\sigma_0+1)}\right] \leq R_0\sup_{k \in X}E_k^{\mu,\nu}\left[R_0^{\sigma_0}\right] \leq R_0 B_0,\ \forall k \neq 0,
\end{eqnarray}
where the equality follows from Proposition 3.4.5(ii) of Meyn and Tweedie \cite{MT1}. Hence we get
\begin{eqnarray} \label{relvalbnd01}
0 < h_1^*[\nu](k) &&= \inf_{\tilde{\mu} \in \SM_1} E_k^{\tilde{\mu},\nu}\left[e^{\theta\sum_{t=0}^{\tau_0}(r_1(X_t,U_t,V_t) - \Lambda_1^*[\nu])}\right] \leq R_0 B_0,\ \forall k \neq 0.
\end{eqnarray}
 Let $\tilde{\mu}[\nu] \in \SM_1$ be any minimizer in (\ref{relval01}). Such a minimizer always exists by (\ref{relvalbnd00}), Proposition \ref{contlawstrat} and the compactness of $\SM_1$. Hence, for all $k \neq 0$, we have using Jensen's inequality and
 part (a) of Proposition \ref{lpstblaprecur}
\begin{eqnarray} \label{relvalbnd02}
 &&h_1^*[\nu](k) = E_k^{\tilde{\mu}[\nu],\nu}\left[e^{\theta\sum_{t=0}^{\tau_0}(r_1(X_t,U_t,V_t) - \Lambda_1^*[\nu])}\right] %\nonumber \\ &&
 = E_k^{\tilde{\mu}[\nu],\nu}\left[e^{\theta\sum_{t=0}^{\sigma_0}(r_1(X_t,U_t,V_t) - \Lambda_1^*[\nu])}\right] \nonumber \\ &&\geq  E_k^{\tilde{\mu}[\nu],\nu}\left[e^{-2\theta||r_1||_\infty (\sigma_0+1)}\right] %\nonumber \\ &&
 \geq  E_k^{\tilde{\mu}[\nu],\nu}\left[e^{-3\Theta||r_1||_\infty (\sigma_0+1)}\right] \geq \frac{1}{R_0}E_k^{\tilde{\mu}[\nu],\nu}\left[\frac{1}{R_0^{\sigma_0}}\right]  \geq \frac{1}{R_0^{E_k^{\tilde{\mu}[\nu],\nu}[\sigma_0] + 1}} \nonumber \\ &&\geq \frac{1}{R_0 E_k^{\tilde{\mu}[\nu],\nu}\left[R_0^{\sigma_0}\right]} %\nonumber \\ &&
 \geq \frac{1}{R_0 \sup_{j \in X}E_j^{\tilde{\mu}[\nu],\nu}\left[R_0^{\sigma_0}\right]} \geq \frac{1}{R_0 B_0}.
\end{eqnarray}
Note that $h_1^*[\nu](0) = \inf_{\xi \in \PR(U)} \int_V \int_U e^{\theta r_1(0,u,v) - \theta \Lambda_1^*[\nu]} \xi(du) \nu\left[dv|0\right]$ obviously obeys the bound. The joint continuity of $E_k^{\mu,\nu}\left[e^{\theta\sum_{t=0}^{\tau_0}(r_1(X_t,U_t,V_t) -\Lambda_1^*[\nu])}\right]$ under $(\mu,\nu)$ can be proved as in the proof of Proposition \ref{eqmingpeergcost}. Hence the continuity of
$h_1^*[\nu]$ in $\nu$ follows from (\ref{relvalbnd00}) and the compactness of $\SM_1$. \hfill $\Box$

Now we state and prove following theorem which is a key step for proving the existence of Nash equilibria.

\begin{theorem} \label{relvalfnrepmpeuniq}
Under \textbf{(A1)} - \textbf{(A3)} and for $(\mu,\nu) \in \SM_1
\times \SM_2$, $(\ln h_1^*[\nu],\hat{\lambda}_1[\nu])$ with $h_1^*[\nu](0) = \int_V \int_U e^{\theta r_1(0,u,v) - \theta \Lambda_1^*[\nu]} w_1^*[\nu](0)(du) \nu\left[dv|0\right]$, where $w_1^*[\nu]$ is defined in (\ref{mpearg01}), is the unique solution in $B(X) \times \R_+$ to (\ref{HJIergupper}). Similarly, $(\ln h_2^*[\mu],\hat{\lambda}_2[\mu])$ with $h_2^*[\mu](0) = \int_U \int_V e^{\theta r_2(0,u,v) - \theta \Lambda_2^*[\mu]} \mu\left[du|0\right] w_2^*[\mu](0)(dv)$, where $w_2^*[\mu]$ is defined in (\ref{mpearg02}), is the unique solution to (\ref{HJIerglower}) again in $B(X) \times \R_+$.
\end{theorem}

\noi \textbf{Proof:} We shall prove this for player I when
player II plays the strategy $\nu$. The other case is proved
analogously. We first prove that $(\ln
h_1^*[\nu],\hat{\lambda}_1[\nu])$ is a solution to
(\ref{HJIergupper}). By Lemma \ref{finrelvalmpeineq}, $h_1^*[\nu]$
as defined in (\ref{relval01}) is finite. Consider $w_1^*[\nu]$
defined in (\ref{mpearg01}) with this $h_1^*[\nu]$.
%Consider the strategy $\tilde{w}_1[\nu] \defn (w_1^*[\nu],w_1^*[\nu],\ldots)$ and let
Let
\begin{equation} \label{mpesol}
f_1^*[\nu](k) \defn
E_k^{w^*_1[\nu],\nu}\left[e^{\theta\sum_{t=0}^{\tau_0}(r_1(X_t,U_t,V_t)
- \Lambda_1[\nu](w^*_1[\nu]))}\right]
\end{equation}
%where $\{w^*_1[\nu]\}$ is generated by
under the stationary strategy $w^*_1[\nu]$.
Note that, for $k \neq 0$, we can get as in (\ref{relvalbnd01}) of Lemma \ref{relvalfnbndlem},
\begin{equation} \label{mpesolbnd}
0 < f_1^*[\nu](k)  \leq R_0 B_0 < \infty.
\end{equation}
It follows from the Markov property that
\begin{eqnarray} \label{markovprobkernz}
&&\int_V \int_U \left(\sum_{j \in X} f_1^*[\nu](j)q(j|k,u,v)\right) w_1^*[\nu](k)(du)\nu\left[dv|k\right] %\nonumber \\ &&
=E_k^{w^*_1[\nu],\nu}\left[e^{\theta\sum_{t=1}^{\sigma_0}(r_1(X_t,U_t,V_t)
- \Lambda_1[\nu](w^*_1[\nu]))}\right] \nonumber \\ && =
\int_V \int_U e^{-\theta r_1(k,u,v) + \theta \Lambda_1[\nu](w^*_1[\nu])} f_1^*[\nu](k) w^*_1[\nu](k)(du) \nu\left[dv|k\right]
,\ k \neq 0.
\end{eqnarray}
Since, $f_1^*[\nu](0) = \int_V \int_U e^{\theta r_1(0,u,v) - \theta
\Lambda_1[\nu](w^*_1[\nu])} w^*_1[\nu](0)(du) \nu\left[dv|0\right]$, it follows that the following
multiplicative Poisson equation holds:
\begin{eqnarray} \label{probkermpe}
&&\int_V \int_U \left(\sum_{j \in X} f_1^*[\nu](j)q(j|k,u,v)\right) w_1^*[\nu](k)(du)\nu\left[dv|k\right] \nonumber \\ =
&&\int_V \int_U e^{-\theta r_1(k,u,v) + \theta \Lambda_1[\nu](w^*_1[\nu])}
f_1^*[\nu](k) w_1^*[\nu](k)(du)\nu\left[dv|k\right], \forall k.
\end{eqnarray}
Since by (\ref{lambdabnd}) in Proposition \ref{eqmingpeergcost}
$\Lambda_1[\nu](\{w^*_1[\nu]\}) < \infty$ and by
(\ref{mpesolbnd}) $0 < f_1^*[\nu](k) < \infty$, under the transition
kernel $q(\cdot|k,u,v)$, we can define the
``twisted kernel" $\breve{q}_1(\cdot|k,w_1^*[\nu](k),\nu(k))$ as
follows:
\begin{equation} \label{twistker}
\breve{q}_1(j|k,w_1^*[\nu](k),\nu(k)) \defn
\int_U \int_V \frac{e^{\theta r_1(k,u,v) - \theta
\Lambda_1[\nu](w^*_1[\nu])}}{f_1^*[\nu](k)}f_1^*[\nu](j)q(j|k,u,v) w_1^*[\nu](k)(du)\nu\left[dv|k\right],\
j \in X.
\end{equation}
Note that the kernel $\breve{q}_1$ is indeed stochastic i.e.,
$\sum_{j \in X} \breve{q}_1(j|k,w_1^*[\nu](k),\nu(k)) = 1,\
\forall k$ since (\ref{probkermpe}) holds. Hence, $\breve{q}_1$ is
the transition kernel for some Markov chain $\{\breve{X}_t\}$.
Since, by (\ref{lambdabnd}), for any $-||r_1||_{\infty} \leq \Lambda <
\Lambda_1[\nu](w^*_1[\nu])$,
\begin{eqnarray*}
%&&
0 < E_0^{w^*_1[\nu],\nu}\left[e^{\theta\sum_{t=0}^{\sigma_0-1}(r_1(X_t,U_t,V_t) - \Lambda)}\right]
\leq E_0^{w^*_1[\nu],\nu}\left[e^{3\Theta||r_1||_\infty \sigma_0}\right]  %\\ &&
\leq \sup_{k \in X}E_k^{w^*_1[\nu],\nu}\left[R_0^{\sigma_0}\right] \leq B_0 < \infty,
\end{eqnarray*}
it follows from Theorem 4.1(i)-(ii) of Balaji and Meyn \cite{BaM1} that the chain
$\{\breve{X}_t\}$ under $\breve{q}_1$ is geometrically recurrent.
Hence it follows from Theorem 4.2(i) of Balaji and Meyn \cite{BaM1} that
$(f_1^*[\nu],\Lambda_1[\nu](\{w^*_1[\nu]\}))$ is a solution to
the multiplicative Poisson equation:
\begin{equation} \label{mpeorig}
e^{\theta\lambda} \psi(k) = \int_U \int_V e^{\theta r_1(k,u,v)}\left(\sum_{j \in X}\psi(j) q(j|k,u,v)\right)w_1^*[\nu](k)(du)\nu\left[dv|k\right].
\end{equation}
It follows from (\ref{mpearg01}), the fact that $h_1^*[\nu] > 0$ by definition and (\ref{mpineq}) in Lemma \ref{finrelvalmpeineq} that, for all $k \in X$,
\begin{eqnarray} \label{verify}
&&\inf_{\xi \in \PR(U)} \int_V \int_U e^{\theta r_1(k,u,v)}\left(\sum_{j\in X}
h_1^*[\nu](j)q(j|k,u,v)\right)\xi(du)\nu\left[dv|k\right]\nonumber \\ &&=
\int_V \int_U e^{\theta r_1(k,u,v)}\left(\sum_{j \in X}h_1^*[\nu](j)
q(j|k,u,v)\right)w_1^*[\nu](k)(du)\nu\left[dv|k\right] \nonumber \\ &&\leq
e^{\theta\Lambda_1^*[\nu]} h_1^*[\nu](k) \leq
e^{\theta\Lambda_1[\nu](w^*_1[\nu])} h_1^*[\nu](k),
\end{eqnarray}
as $\Lambda_1^*[\nu] \leq \Lambda_1[\nu](w^*_1[\nu])$ (see
(\ref{generalgpe01}) and (\ref{mingpe01})). Hence, by Theorem 4.2(i)-(ii) of Balaji and Meyn
\cite{BaM1}, we obtain
\begin{equation} \label{mpesoleq}
\frac{h_1^*[\nu](k)}{h_1^*[\nu](0)} = \frac{f_1^*[\nu](k)}{f_1^*[\nu](0)},
\end{equation}
and
\begin{eqnarray} \label{finalverify}
&&\inf_{\xi \in \PR(U)} \int_U \int_V e^{\theta r_1(k,u,v)}\left(\sum_{j\in X}
h_1^*[\nu](j)q(j|k,u,v)\right)\xi(du)\nu\left[dv|k\right] \nonumber \\ &&= e^{\theta\Lambda_1^*[\nu]}
h_1^*[\nu](k) = e^{\theta\Lambda_1[\nu](w^*_1[\nu])}
h_1^*[\nu](k),\ \forall k \in X.
\end{eqnarray}
Hence $\Lambda_1^*[\nu] = \Lambda_1[\nu](w^*_1[\nu])$.
Since, by Lemma \ref{relvalfnbndlem}, $0 < \sup_{k \in X} |\ln
h_1^*[\nu](k)| \leq \ln(R_0 B_0) < \infty$ as by definition both $R_0, B_0 > 1$, uniqueness follows
from Corollary 2.3 of Di Masi and Stettner \cite{DS1} and the condition $h_1^*[\nu](0)
= f_1^*[\nu](0) = \inf_{\xi \in \PR(U)}\int_U \int_V e^{\theta r_1(0,u,v) - \theta \Lambda_1^*[\nu]} \xi(du) \nu\left[dv|0\right] $. It
follows from Proposition 1.1 of Di Masi and Stettner \cite{DS1} (see also Theorem 2.1 of Hern\'{a}ndez-Hern\'{a}ndez and Marcus \cite{HM1}) that
\begin{equation} \label{eqmingpestochrep01}
\Lambda_1^*[\nu] = \hat{\lambda}_1[\nu]
\end{equation}
and hence,
\begin{equation} \label{relvalstochrep01}
h_1^*[\nu](k) =
E_k^{\{w^*_1[\nu]\},\nu}\left[e^{\theta\sum_{t=0}^{\tau_0}(r_1(X_t,U_t,V_t)
- \hat{\lambda}_1[\nu])}\right], \forall k \in X.
\end{equation}
Similarly, it can be shown for player II that
\begin{equation} \label{eqmingpestochrep01}
\Lambda_2^*[\mu] = \hat{\lambda}_2[\mu]
\end{equation}
and
\begin{equation} \label{relvalstochrep02}
h_2^*[\mu](k) =
E_k^{\mu,\{w^*_2[\mu]\}}\left[e^{\theta\sum_{t=0}^{\tau_0}(r_2(X_t,U_t,V_t)
- \hat{\lambda}_2[\mu])}\right], \forall k \in X,
\end{equation}
  under the stationary strategy $w_2^*[\mu]$ as defined in (\ref{mpearg02}) with $h_2^*[\mu]$ as defined in (\ref{relval02}). \hfill $\Box$

Hence, by (\ref{discminsel1erg}) and (\ref{discminsel2erg}), for
 any $(\mu,\nu) \in \SM_1 \times \SM_2$,
there exists measurable (minimizing) selectors (see Bene\v{s} \cite{Be1})
\[(\hat{\mu}[\nu],\hat{\nu}[\mu]) : X \mapsto \PR(U) \times \PR(V)\] satisfying
\begin{eqnarray} \label{discminsel1erg2}
&&\inf_{\xi \in \PR(U)} \int_V \int_U e^{\theta r_1(k,u,v)}\left(\sum_{j \in X}e^{\ln h^*_1[\nu](j)} q(j|k,u,v)\right)\xi(du)\nu\left[dv|k\right] \nonumber \\ &&=
\int_V \int_U e^{\theta r_1(k,u,v)} \left(\sum_{j \in X} e^{\ln h^*_1[\nu](j)}q(j|k,u,v)\right)\hat{\mu}[\nu](k)(du)\nu\left[dv|k\right],
\end{eqnarray}
implying, by Theorem \ref{relvalfnrepmpeuniq} and (\ref{HJIergupper}),
\begin{equation} \label{HJIerguppersol}
e^{\theta\hat{\lambda}_1[\nu] + \ln h_1^*[\nu](k)} = \int_V \int_U e^{\theta
r_1(k,u,v)} \left(\sum_{j \in X} e^{\ln
h^*_1[\nu](j)}q(j|k,u,v)\right)\hat{\mu}[\nu](k)(du)\nu\left[dv|k\right],
\end{equation}
and similarly
\begin{eqnarray} \label{discminsel2erg2}
&&\inf_{\chi \in \PR(V)} \int_U \int_V e^{\theta r_2(k,u,v)} \left(\sum_{j \in X} e^{\ln h^*_2[\mu](j)} q(j|k,u,v)\right)\mu\left[du|k\right]\chi(dv)\nonumber \\
&&= \int_U \int_V e^{\theta r_2(k,u,v)} \left(\sum_{j \in X} e^{\ln h^*_2[\mu](j)}\right)q(j|k,u,v)\mu\left[du|k\right]\hat{\nu}[\mu](k)(dv),
\end{eqnarray}
implying, by Theorem \ref{relvalfnrepmpeuniq} and (\ref{HJIerglower}),
\begin{equation} \label{HJIerglowersol}
e^{\theta\hat{\lambda}_2[\mu] + \ln h_2^*[\mu](k)} = \int_U \int_V e^{\theta r_2(k,u,v)} \left(\sum_{j \in X} e^{\ln h^*_2[\mu](j)}q(j|k,u,v)\right)\mu\left[du|k\right]\hat{\nu}[\mu](k)(dv).
\end{equation}

Using the minimizing selectors $\hat{\mu}[\nu]$ and $\hat{\nu}[\mu]$ in (\ref{discminsel1erg2}) and (\ref{discminsel2erg2}) as our optimal responses for given $(\mu,\nu) \in \SM_1 \times \SM_2$, %(\ref{minseq01erg}) and (\ref{minseq02erg}) resp.,
we now prove the existence of Nash equilibria in the class of stationary Markov strategies for the ergodic cost criterion (\ref{riskergpayoff}).
%we now state and prove the result on the existence of a fixed point of the map $\HL$ defined in (\ref{optresponsedisc}).

%\begin{proposition} \label{Nashlem02}
%Given $\theta \in (0,\Theta]$, the map $\HL$ has a fixed point $(\tilde{\mu}(\theta,\cdot),\tilde{\nu}(\theta,\cdot))$.
%\end{proposition}

\begin{theorem} \label{existNasherg}
Given $\theta \in (0,\Theta]$, there exists a pair of Nash
equilibrium strategy in $\SM_1 \times \SM_2$ for the game
corresponding to the cost criterion (\ref{riskergpayoff}).
\end{theorem}

%\noi \textbf{Proof:} Using (\ref{riskergpayoff}), (\ref{stochrep01erg}) and (\ref{stochrep02erg}), this result follows directly from Proposition \ref{Nashlem02} above. \hfill $\Box$

\noi \textbf{Proof:} Given $(\mu,\nu) \in \SM_1 \times \SM_2$,
$\HL$ as defined in \eqref{optresponseerg} can be proved to be
a map with non-empty, closed and convex values as in Theorem \ref{existexpvalue}.
%both $\HL_1(\nu)$ and $\HL_2(\mu)$ are non-empty and convex and hence $\HL(\mu,\nu)$ is also non-empty and convex.The closure property of $\HL(\mu,\nu)$ follows from the joint continuity of $\hat{q}(\cdot|\cdot,\mu,\nu)$ under $(\mu,\nu)$. Hence,
%Now, by Proposition \ref{contlawstrat}, the law of the process $\{X_t\}$ driven by (stationary) Markov strategies is jointly continuous in these strategies.
Now, by Proposition \ref{contlawstrat}, $\LA^{\mu,\nu}\{X_t, t \geq 0\}$
is jointly continuous in $(\mu,\nu) \in \SM_1 \times \SM_2$.
Note that $\hat{\lambda}_1[\nu]$ as obtained in
(\ref{stochrep01erg}) has a minimizing stationary strategy
(\ref{discminsel1erg}) and $\hat{\lambda}_2[\mu]$ as obtained in
(\ref{stochrep02erg}) has a minimizing stationary strategy
(\ref{discminsel2erg}), i.e., the minimizations in (\ref{stochrep01erg}) and (\ref{stochrep02erg})
can be effectively considered over $\SM_1$ and $\SM_2$ (resp.) which are compact sets. Also, $r_i,\ i=1,2$ are bounded functions.
Hence, %by Proposition \ref{contlawstrat} and Skorohod's theorem (see Theorem 2.2.2 of Borkar \cite{B1}),
by Berge Maximum Theorem (see Theorem 17.31 of Aliprantis and Border \cite{AliBor1}), $\hat{\lambda}_1[\nu]$ and
$\hat{\lambda}_2[\mu]$ are continuous in $\nu$ and $\mu$
respectively. Similarly, the continuity of $h_1^*[\nu]$ (resp.\ $h_2^*[\mu]$) follows from (\ref{relvalstochrep01}) (resp.\ (\ref{relvalstochrep02})) and Lemma \ref{relvalfnbndlem}. %and Proposition \ref{contlawstrat}.
Now consider a sequence $\{(\mu_k,\nu_k)\}_{k \in \mathbb{N}}$ in
$\SM_1 \times \SM_2$ converging to $(\mu_\infty,\nu_\infty) \in
\SM_1 \times \SM_2$. Let $\HL(\mu_k,\nu_k) \ni
(\tilde{\mu}_k,\tilde{\nu}_k) \ra
(\tilde{\mu}_\infty,\tilde{\nu}_\infty) \in \SM_1 \times \SM_2$.
From the continuity of
$q(\cdot|\cdot,u,v),(h_1^*[\nu],\hat{\lambda}_1[\nu])$
and $(h_2^*[\mu],\hat{\lambda}_2[\mu])$ in $\mu,\nu$, it follows
that $(\tilde{\mu}_\infty,\tilde{\nu}_\infty) \in
\HL(\mu_\infty,\nu_\infty)$. Hence the map $\HL$ is upper
semi-continuous. Then, by Theorem 1 of Fan \cite{F1}, there exists a
fixed point $(\mu^*,\nu^*) : X \mapsto \PR(U) \times \PR(V)$ of the
map $\HL$ defined as
\begin{equation} \label{expergcostfp}
(\mu^*(k),\nu^*(k)) \defn
(\hat{\mu}[\nu^*](k),\hat{\nu}[\mu^*](k))
\end{equation}
where $\hat{\mu}[\nu^*]$ (resp.\ $\hat{\nu}[\mu^*]$) is any
minimizing selector in (\ref{discminsel1erg2}) (resp.
(\ref{discminsel2erg2})). Then $(\mu^*,\nu^*) \in \SM_1 \times
\SM_2$ is a Nash equilibrium strategy for the cost criterion
(\ref{riskergpayoff}), i.e., for all $k \in X$
\[\hat{\lambda}_1[\nu^*] = \beta^{\mu^*,\nu^*}_1(\theta,k) \leq \beta^{\mu,\nu^*}_1(\theta,k),\ \forall \mu \in \Pi_1, \]
and
\[\hat{\lambda}_2[\mu^*] = \beta^{\mu^*,\nu^*}_2(\theta,k) \leq \beta^{\mu^*,\nu}_2(\theta,k),\ \forall \nu \in \Pi_2. \]
\hfill $\Box$

\begin{comment}
We now provide an illustrative example.

\begin{example} \label{nashergexmpl}
\end{example}
\end{comment}

\section{Conclusions.} \label{conclude}
We have established Nash equilibria for nonzero-sum discrete-time stochastic
games on a countable state space with risk-sensitive cost criteria. We have studied both
discounted and ergodic cost criteria on the infinite time horizon.
For discounted cost criterion, we have established Nash
equilibria under fairly general conditions. %when the cost is
%independent of the controls of the players \textbf{(A0)}.
For ergodic cost criterion, however, we have established our results
under additional criteria of irreducibility \textbf{(A1)},
stability \textbf{(A2)} and smallness of cost \textbf{(A3)} all of which
are standard in the corresponding literature.
A natural continuation of this work is to analyze such discrete-time
games on general state spaces. Another direction is to investigate
such games in continuous time in the set-up of stochastic differential
games.

% Appendix here
% Options are (1) APPENDIX (with or without general title) or
%             (2) APPENDICES (if it has more than one unrelated sections)
% Outcomment the appropriate case if necessary
%
%\begin{APPENDIX}
%\end{APPENDIX}
%
%   or
%
% \begin{APPENDICES}
% \section{<Title of Section A>}
% \section{<Title of Section B>}
% etc
% \end{APPENDICES}

\section*{Acknowledgments.}
%The first author gratefully acknowledges the support of this research in part
 %by Department of Science and Technology, Govt.\ of India
 %research grant no: SR/FTP/MS-08/2009 and in part by Research Chair grant no: 190
 %of IIM Bangalore.
 The work of the second author is supported in part by a grant from
 UGC Centre for Advanced Study.

\bibliographystyle{nonumber}

\end{document}